\documentclass[final,reqno]{amsart}
\usepackage{pdfsync}
\usepackage{amsmath,amssymb}
\usepackage{enumerate}
\usepackage{xspace}
\usepackage{array}
\usepackage{tabularx}
\usepackage{booktabs}
\usepackage{subfigure}
\usepackage{caption}
\usepackage{fullpage}
\usepackage{cite}
\usepackage{caption}
\usepackage{james}
\usepackage{cleveref}
\graphicspath{{figures/}}

\newlength\mylen
\newcommand\myinput[1]{%
  \settowidth\mylen{\KwInput{}}%
  \setlength\hangindent{\mylen}%
  \hspace*{\mylen}#1}
\newcommand{\KwAnd}{\textup{and}~}

\let\oldnl\nl
\newcommand{\nonl}{\renewcommand{\nl}{\let\nl\oldnl}}

\hypersetup{colorlinks, linkcolor=blue, citecolor=blue, urlcolor=blue}

\definecolor{BLUE}{RGB}{0,0,255} 
\definecolor{RED}{RGB}{255,0,0}

\newcommand{\revised}[1]{{#1}}

\usepackage{mathtools}
\mathtoolsset{showonlyrefs=true} 

\numberwithin{equation}{section}

\author{ James Jackaman }
\address{ James Jackaman \thanks{ Department
    of Mathematical Sciences,
    NTNU, 7491 Trondheim, Norway {\tt{james.jackaman@ntnu.no}}.  }}

\author{Scott MacLachlan}
\address{
  Scott MacLachlan
  \thanks{
    Department of Mathematics and Statistics, Memorial University of
    Newfoundland, St.\ John's, NL, A1C 5S7, Canada
    {\tt{smaclachlan@mun.ca}}.
  }}

\thanks{This work was partially supported by an NSERC Discovery Grant
  (SM and JJ) and the ERCIM Alain Bensoussan Fellowship Programme
  (JJ)}

\title[]%
{Constraint-Satisfying Krylov Solvers for Structure-Preserving
  Discretisations} \date{\today}

\begin{document}

\maketitle

\begin{abstract}
  A key consideration in the development of numerical schemes for
  time-dependent partial differential equations (PDEs) is the ability
  to preserve certain properties of the continuum solution, such as
  associated conservation laws or other geometric structures of the
  solution.  There is a long history of the development and analysis
  of such \textit{structure-preserving} discretisation schemes,
  including both proofs that standard schemes have
  structure-preserving properties and proposals for novel schemes that
  achieve both high-order accuracy and exact preservation of certain
  properties of the continuum differential equation.  When coupled
  with implicit time-stepping methods, a major downside
  to these schemes is that their structure-preserving properties
  generally rely on \textit{exact} solution of the (possibly
  nonlinear) systems of equations defining each time-step in the
  discrete scheme.  For small systems, this is often possible (up to
  the accuracy of floating-point arithmetic), but it becomes
  impractical for the large linear systems that arise when considering
  typical discretisations of space-time PDEs.  In this paper, we
  propose a modification to the standard flexible generalised minimum
  residual (FGMRES) iteration that enforces selected constraints on
  approximate numerical solutions.  We demonstrate its application to
  both systems of conservation laws and dissipative systems.
\end{abstract}


\section{Introduction} \label{sec:introduction}

Since the invention of modern digital computers, there has been
intense study of numerical methods for the approximation of solutions
to systems of partial differential equations (see, for example,
\cite{FHHarlow_JEWelch_1965a}).  Within this study has arisen a
specialisation of methods for continuum PDEs that have associated
geometric properties, such as Hamiltonian and reversible systems,
problems posed over Lie groups, or systems with other structures
deemed noteworthy.  This field is quite mature, with numerous research
papers and textbooks published in the area known as \textit{geometric
  numerical integration}~\cite{Sanz-SernaCalvo:1994,
  LeimkuhlerReich:2004, HairerLubichWanner:2006, BlanesCasas:2016},
but also with a continuing interest in the development and analysis of
new methods within the class continuing to the present (cf., for
example, \cite{Holec:2022, SatoMiyatakeButcher:2022,
  Bogfjellmo:2022}). The present work is motivated by the
impracticality of an underlying assumption in this work, that the
linear and nonlinear systems of equations associated with these
discretisations can be effectively solved to high-enough accuracy to
guarantee preservation of geometric structure (at least up to machine
precision of floating-point arithmetic).

The solution of systems of equations is generally required by any
implicit numerical scheme for time-dependent PDEs.  While some
geometric schemes can utilise explicit integrators
(e.g., \cite{Tao:2016}), schemes typically rely on implicit
discretisation to achieve both stability and their geometric
properties.  When considering numerical solution of small systems of
ordinary differential equations, we rarely account for the true
computational burden of solving linear and nonlinear systems exactly,
namely using LU factorisation to solve the linear systems, and
Newton's method to high tolerance for nonlinear systems.  For systems
of ODEs that arise from the method-of-lines (or other) discretisation
of space-time PDEs, however, these assumptions are often impractical,
due simply to the computational cost of exact solutions of linear
systems, or of solving either linear or nonlinear systems to high
tolerances.

Because of this essential issue, some attention has been drawn to the
development of Krylov methods for the solution of systems
with Hamiltonian or symplectic structure~\cite{Benner:1997,
  Benner:2001, Watkins:2004}.  The inherent idea in this work (and its
generalisations~\cite{Benner:2011, Li:2019}) is to build
the Krylov space so that all vectors in the Krylov space satisfy the
expected properties of the solution.  In the Hamiltonian case, we
consider a symmetric matrix, $H$, of dimension $2m\times 2m$, and
skew-symmetric matrix, $J = \left[\begin{smallmatrix} 0 & I \\ -I &
  0 \end{smallmatrix}\right]$, where the blocks of $J$ are of size $m\times
m$, and the symplectic Lanczos algorithm~\cite{Benner:2011} generates
a $J$-orthogonal basis for the Krylov spaces of even dimension.
In~\cite{Li:2019}, it is shown that this can be used to generate
energy-preserving approximate numerical solutions of the ODE system
$\frac{d\vec{y}}{dt} = JH\vec{y}$.  To our knowledge, these approaches
have not been extended beyond variations on Hamiltonian/symplectic
structure, or to broader classes of conservation laws.  In this paper,
we propose an alternate approach that is more costly, but also more
flexible, allowing us to work with arbitrary conservation laws,
regardless of their origin or geometric structure.

The essential ingredients of our approach are as follows.  We consider
classes of space-time PDEs that have associated geometric properties,
either in terms of related conserved quantities or dissipation laws.
For any given problem, we rely upon a structure-preserving
discretisation; for our problems, this will take the form of a spatial
semi-discretisation that preserves the continuum properties (see, for
example, \cite{CotterShipton:2012, self:thesis, Holec:2022}) coupled
with an implicit-in-time Runge-Kutta (RK) method that preserves these
properties now in the discrete setting (at each timestep)
(cf.~\cite{SatoMiyatakeButcher:2022}).  At the discrete level, this
gives us a nonlinear or linear system to be solved at each timestep
(we consider only the linear case herein), defining the approximate
solution at the next timestep in terms of that at the current
timestep, along with a list of constraint functionals (e.g., linear or bilinear
forms, although more general constraints are possible) that have
prescribed values for the approximate solution at the next time step.
In this paper, we describe a modification to the standard flexible
GMRES (FGMRES) algorithm~\cite{SaadSchultz:1986, YSaad_2003a} that
provides an approximate solution to the linear system that
simultaneously satisfies (or attempts to satisfy) the constraints
associated with the list of functionals given.

While two variants are described, the main approach here can be
implemented in a relatively non-intrusive way to an existing FGMRES
implementation, requiring only access to the basis for
the (preconditioned) Krylov space once the FGMRES convergence criterion is (close to)
satisfied.  Notably, the method can be implemented without
restrictions on the use of any (sensible) preconditioning strategy.
Furthermore, the approach is viable for any Krylov method that
maintains (i.e., stores) a basis for the Krylov space, although it is
particularly well-suited to the structure of calculation within GMRES
and FGMRES.  The main expense in the method is the solution of a
constrained minimisation problem, for which we make use of \revised{a
sequential quadratic programming approach~\cite{DKraft_1988}} as implemented in
SciPy~\cite{2020SciPy-NMeth}.  We note that both our intent and
approach is quite different from that in~\cite{BirkenLinders:2021, LindersBirken2022},
which focus on local and global conservation (of mass) for hyperbolic conservation laws by \textit{unmodified} iterative
methods for solution of the resulting linear and nonlinear systems.  \revised{We also note the similarity between our approach and that used in~\cite{Frommer:1998}, where solutions to a given linear system and its shifted counterpart are computed from the same Krylov space by retaining similar information with additional computation.}

The remainder of this paper is organised as follows.
In Section~\ref{sec:background}, we give an overview of structure-preserving
discretisation based on the finite-element plus implicit RK
methodology.  Then, in Section~\ref{sec:methodology}, we present the main
contribution here, of the modified constraint-satisfying form of the
FGMRES algorithm, and some discussion of its limitations and possible
extensions.  Numerical examples are presented in Section~\ref{sec:examples},
followed by concluding remarks in Section~\ref{sec:conclusion}.

\section{Structure-Preserving Discretisation}\label{sec:background}

Many PDEs (or systems of PDEs) possess some underlying structure, be
it symplectic, variational, symmetry preserving or conserved
quantities, and much work has gone into developing discretisation
frameworks that preserve such
structures~\cite{HairerLubichWanner:2006, MarsdenWest:2001,
  BlanesCasas:2016, Mansfield:2013}. Here we focus on the case of
linear (systems of) PDEs, where the structure-preserving properties of
such discretisation depends on suitable solution of the underlying
linear systems in an implicit time discretisation strategy.  We note
that the ideas outlined here may be generalised to nonlinear systems
of PDEs and nonlinear solvers, as discussed in
Remark~\ref{rem:nonlinear}.

For simplicity, we focus on numerical schemes that
possess associated conservation laws, but note that the methodology
can be applied to any underlying structure that can be expressed in a
similar framework (including the case of
a dissipative system discussed in Section~\ref{sec:background:heat}).  We
consider discretisation using the method-of-lines approach, coupling
well-chosen spatial finite-element methods with implicit RK
temporal integrators.  The specific choice of the spatial
finite-element method depends strongly on the problem at hand;
it is rare that standard choices of continuous piecewise polynomial
(Lagrange) spaces yield \revised{Galerkin} discretisations with significant conservative
properties\revised{, although summation-by-parts methods provide an alternative that does just this, cf.~\cite{doi:10.1137/120890144,GASSNER201639}}. In general, we consider a continuum PDE or system of
PDEs, $u_t + Lu = g$, posed on the time interval $[0,T]$ and spatial
domain, $\Omega$, with suitable boundary conditions on
$\partial\Omega$.  We pass from the PDE form to a variational
formulation by multiplying by a (scalar or vector) test function,
$\phi$, and (typically) integrating by parts on some or all of the
terms in $\langle Lu,\phi \rangle$.  Then, we choose suitable discrete
subspace(s) for $u$ and $\phi$, leading to the semi-discretised form
of the PDE, which we write as $\vec{z}_t = \revised{\vec{h}}(t,\vec{z})$ for
$\vec{z}(t)\in \mathbb{R}^d$ and $t \in [0,T]$.  Details of this
process are described for three model problems in the following
subsections.  We note that there is a one-to-one correspondence
between the discrete values of $\vec{z}(t)$ at time $t$ and a function
on $\Omega$ that is the finite-element approximation to $u(t,x)$ for
$x\in\Omega$.  We also assume that this semi-discretisation process
preserves some number of continuum conservation laws that, after
discretisation, can be expressed as $\frac{d}{dt} g_\ell(\vec{z}(t)) =
0$.

For the temporal discretisation, we consider the symplectic
\revised{Gauss-Legendre} family of implicit RK methods (GLRK).  For a problem of the form $\vec{z}_t = \revised{\vec{h}}(t,\vec{z})$, we write
the $s$-stage symplectic RK method as
\begin{equation} \label{eqn:symrk}
    \vec{z}^{n+1} = \vec{z}^n + \dt{n} \sum_{i=1}^s b_i \vec{k}^i \text{ with }
    \vec{k}^i = \revised{\vec{h}}\bc{t^n+c_i\dt{n}, \vec{z}^n + \dt{n} \sum_{j=1}^s a_{ij} \vec{k}^j}
    ,
\end{equation}
where $\vec{z}^n$ is the approximation to $\vec{z}(t_n)$ for $t_n\in
[0,T]$, $\dt{n} = t_{n+1}-t_n$, and the values $\{a_{ij}\}$,
$\{b_i\}$, and $\{c_i\}$ come from the \textit{Butcher Tableau} for
the RK scheme.  The vectors $\vec{k}^i$ denote an approximation to the
derivative of $\vec{z}(t)$ at the stage times $t^n+c_i\dt{n}$.  Here,
we take $\{c_i\}$ to be the zeroes of the shifted Legendre
polynomial $\frac{d^s}{dx^s}\bc{x^s(x-1)^s}$, and the other coefficients of
\eqref{eqn:symrk} are described by
\begin{equation}
  a_{ij} = \int_0^{c_i} \mathcal{L}_j(\tau) \di{\tau},
  \qquad
  b_i = \int_0^1 \mathcal{L}_i(\tau) \di{\tau}
,
\end{equation}
where $\mathcal{L}_i$ is the $i$-th Lagrange polynomial. We note that,
as proven in \cite[IV.2.1]{HairerLubichWanner:2006}, this family of
methods preserves all quadratic invariants, meaning that
$g_\ell(\vec{z}^{n+1}) = g_\ell(\vec{z}^n)$ if $g_\ell(\vec{z})$ is a
polynomial in $\vec{z}$ of degree at most two. For simplicity, we typically employ the
lowest-order method in this family, the implicit midpoint method. We
note that this, in turn, is equivalent to the popular \CN temporal
discretisation methodology, under the assumption that
$\revised{\vec{h}}(t,\vec{z})$ is a linear operator.

Practically, \eqref{eqn:symrk} results in a coupled system of
equations to solve for $\{\vec{k}^i\}$ at each time step.  While, in
general, this system is nonlinear (due to dependence on the nonlinear
function $\revised{\vec{h}}(t,\vec{z})$), we consider the case of linear PDEs,
leading to linear systems of equations of dimension $\mathbb{R}^{sd}$,
with one vector $\vec{k}^i \in \mathbb{R}^d$ for each stage in the RK
scheme.  Using standard notation, we write this linear system as
$\revised{\A}\vec{x}=\revised{\vec{f}}$, suppressing dependence on the time step and noting
that matrix $\revised{\A} \in \mathbb{R}^{sd\times sd}$ is not the same as the
coefficients $\{a_{ij}\}$ in the Butcher Tableau.  We note that in the
special case of the \CN discretisation, this linear system can also be
solved directly for $\vec{z}^{n+1}$ instead of the single stage value
$\vec{k}^1$, which we do in that case.

\subsection{The Linear Korteweg-De Vries Equation} \label{sec:background:lkdv}

The linear Korteweg-De Vries (KdV) equation is given by
\begin{equation} \label{eqn:lkdv}
  u_t + u_x  + u_{xxx} = 0,
\end{equation}
where we consider $t\in [0,T]$ and $x\in \Omega = [0,X]$, with $u(t,x)$ being
periodic with period $X$, and initial condition $u(0,x)=u_0(x)$.  Just
as in the case of the nonlinear KdV equation, this system has several
nontrivial invariants.  For the sake of exposition, we focus on
three of these, \revised{namely mass conservation,
$
\ddt \int_\Omega u \di{x} = 0
$, 
momentum conservation,
$
\ddt \int_\Omega \frac12 u^2 \di{x} = 0
$, 
and energy conservation,
$
\ddt \int_\Omega \frac12 u_x^2 - \frac12 u^2 \di{x} = 0
$.
}

Given the dependence on a high-order spatial derivative in
\eqref{eqn:lkdv}, the spatial discretisation requires either 
higher-order finite-element methods with enhanced continuity or
rewriting the equations \revised{to allow} use of standard Lagrange
basis functions.  We choose the latter, rewriting the linear KdV
equation as the equivalent system
\revised{
\begin{equation}
  \begin{split}
  u_t + v_x = 0, \qquad
  v -u - w_x = 0, \qquad
  w -u_x = 0.
  \end{split}
\end{equation}
}%
Partitioning the interval, $\Omega$, as
$0=x_0<\cdots<x_m<\cdots<x_{\M}=X$ with an arbitrary element
$\selement{m} := (x_m, x_{m+1})$, we
define the spatial finite-element space as follows.

\begin{definition}[One dimensional finite-element spaces] \label{def:sfes}
  \index{Spatial finite element spaces}

  Let $\bpoly{q}(\selement{m})$ denote the space of polynomials of
  degree no more than $q$ on the element $\selement{m}$. Then, the
  discontinuous finite-element space $\dpoly{q}$ is given by
  \begin{equation} \label{eqn:dpoly}
    \dpoly{q}
    =
    \{U : \Omega \to \mathbb{R} : \left. U \right|_{\selement{m}} \in \bpoly{q}
    (\selement{m})
    \textrm{ for } m=0,...,{\M}-1 \}
    .
  \end{equation}
  Note that the dimension of $\dpoly{q}$ is
  given by $\bc{q+1}M_x$
\end{definition}
Before introducing the finite element method, we must first define the
\revised{weak} first spatial derivative for the discontinuous function space.
\begin{definition}[Spatial derivative operator] \label{def:gfunc} Let
  $U \in \dpoly{q}$, then we define $\gfunc{} : \dpoly{q} \to \dpoly{q}$ such
  that
  \begin{equation} \label{eqn:gfunc}
    \int_\Omega \gfunc{U} \phi
    =
    \sum_{m=0}^{\M-1} \int_{\selement{m}} U_x \phi 
    - \jump{U_m} \avg{\phi_m}
    \qquad
    \forall \phi \in \dpoly{q}
    ,
  \end{equation}
  where
  \begin{equation}
    \jump{U_m} = \lim_{y \nearrow x_m} U(x) - \lim_{y \searrow x_m}
    U(x)
  \end{equation}
  is the jump in $U(x)$ at $x=x_m$ and
  \begin{equation}
    \avg{U_m} = \frac12 \bc{ \lim_{y \nearrow x_m} U(x) + \lim_{y
        \searrow x_m} U(x)}
  \end{equation}
  the average of the two one-sided limits at
  $x_m$. 
\end{definition}
With this definition in mind, our conservative spatial finite-element
discretisation of \eqref{eqn:lkdv} is given by seeking
$U,V,W \in \dpoly{q}$ such that
\begin{equation} \label{eqn:spatiallkdv}
  \begin{split}
    \int_\Omega \bc{U_t + \gfunc{V}} \phi \di{x}  &= 0
    \qquad \forall \phi \in \dpoly{q} 
    \\
    \int_\Omega \bc{ V - U - \gfunc{W}}  \psi \di{x}  & = 0
    \qquad \forall \psi \in \dpoly{q}
    \\
    \int_\Omega \bc{ W - \gfunc{U}} \chi \di{x} & = 0 
    \qquad \forall \chi \in \dpoly{q}
    ,
  \end{split}
\end{equation}
subject to \revised{initialising $U$, $V$, and $W$ based on some initial data $U(0,x)$}. We note there are
similarities between the first-order system form of the linear KdV
equation and the design of local discontinuous Galerkin methods; see,
for example, \cite{YanShu:2002, HuffordXing:2014}; however, we
use only the central fluxes defined above (see
\cite[\S4]{self:thesis}).  We also note that, in contrast to the
earlier discussion, \eqref{eqn:spatiallkdv} takes the form of a
differential-algebraic equation (DAE), for which Runge-Kutta methods
remain well-defined~\cite{wanner1996solving}.

We fully discretise the system in~\eqref{eqn:spatiallkdv} with a
method of lines approach, partitioning the temporal interval as
$0=t_0<\cdots<t_n<\cdots<t_{\N}=T$, and evaluate the spatial
finite-element discretisation at the discrete points in time $t_n$,
with the initial condition given by a spatial finite-element
formulation at $t_0$. We begin by considering a \CN-type
discretisation which, for linear problems, is equivalent to the
lowest-order symplectic RK method, preserving quadratic invariants.

\begin{definition}[Finite-element discretization for linear KdV] \label{def:lkdv}

  Let $U^j,W^j \in \dpoly{q}$ be given for $j=0,..,n$. Then, we seek
  $U^{n+1},V,W^{n+1} \in \dpoly{q}$ such that
  \begin{equation} \label{eqn:fulllkdv}
    \begin{split}
      \int_\Omega \bc{\frac{U^{n+1}-U^n}{\dt{n}} + \gfunc{V} } \phi \di{x}  &= 0
      \qquad \forall \phi \in \dpoly{q} 
      \\
      \int_\Omega \bc{ V - U^{n+\frac12} - \gfunc{W^{n+\frac12}} } \psi \di{x} & = 0
      \qquad \forall \psi \in \dpoly{q}
      \\
      \int_\Omega \bc{ W^{n+1} - \gfunc{U^{n+1}}  } \chi \di{x} & = 0 
      \qquad \forall \chi \in \dpoly{q}
      ,
    \end{split}
  \end{equation}
  where $U^{n+\frac12} := \frac12 \bc{U^{n+1}+U^n}$ and
  $W^{n+\frac12} := \frac12 \bc{W^{n+1}+W^n}$ for $n =
  0,..,{\N}-1$. We define the initial data $U^0 = \Pi u_0(x)$ where
  $\Pi$ denotes the $L_2$ projection into the finite-element space,
  $\dpoly{q}$, and we initialise $W$ such that $W^0 = \gfunc{U^0}$
  (as required for the conservation properties of the scheme
  \cite[\S4]{self:thesis}). Note that only the intermediate value,
  $V^{n+\frac12}$, is required for the scheme; while we could compute
  this from $V^n$ and $V^{n+1}$, this does not change the
  conservative properties of the discretisation, so we neglect the
  time dependence in this term, noting that it arises implicitly
  through the
  dependence on $U$ and $W$. In addition, we evaluate the third
  equation at $t_{n+1}$ as, for momentum conservation, it is important
  to be able to reformulate a discrete difference of form
  \begin{equation}
    \int_\Omega \bc{ \frac{W^{n+1}-W^n}{\dt{n}}
      - \gfunc{\frac{U^{n+1}-U^n}{\dt{n}}}  } \chi \di{x}
    = 0
    .
  \end{equation}
\end{definition}

This finite element discretisation preserves fully discrete versions
of mass, momentum and energy (see \cite[\S4]{self:thesis}), which may
be expressed as
\begin{equation} \label{eqn:lkdv:laws}
  \begin{split}
    \int_\Omega U^{n+1} \di{x}
    & =
    \int_\Omega U^n \di{x}
    \\
    \int_\Omega \frac12 \bs{U^{n+1}}^2 \di{x}
    & =
    \int_\Omega \frac12 \bs{U^n}^2 \di{x} \\
    \frac12 \int_\Omega \bs{W^{n+1}}^2 - \bs{U^{n+1}}^2 \di{x}
    & =
    \frac12 \int_\Omega \bs{W^n}^2 - \bs{U^n}^2 \di{x}
    .
  \end{split}
\end{equation}
We write the complete finite-element solution at $t=t_{n}$ as
$\vec{Z}^{n} = [U^{n}; V; W^{n}] \in \dpoly{q}^3$ through concatenating the
solution. In turn, this allows us to define the discrete solution vector
$\vec{z}^{n} \in \mathbb{R}^{3(q+1)M_x}$ as the basis coefficients of
$\vec{Z}^{n}$. We can then
express these constraints by defining
$g_1(\vec{z}^n) = \int_\Omega U^n \di{x}$,
$g_2(\vec{z}^n) = \int_\Omega \frac12 \bs{U^n}^2\di{x}$, and
$g_3(\vec{z}^n) = \frac12 \int_\Omega \bs{W^n}^2 - \bs{U^n}^2 \di{x}.$
In order to explicitly write the constraints in terms of the vector
$\vec{z}^{n+1}$, we define the vector $\vec{\omega}\bc{u}\in \mathbb{R}^{3(q+1)M_x}$ such that
\begin{equation}
  \omega(u)_i
  :=
  \int_\Omega \basis{u}_i \di{x}
  ,
\end{equation}
where $\basis{u}_i$ represents the $i$-th basis function for variable $U$ in the trial
space for $\vec{Z} \in\dpoly{q}^3$. Similarly, we define the mass matrices of the vector components
\begin{equation}
  M(u)_{ij}
  :=
  \int_\Omega \basis{u}_i \basis{u}_j \di{x}
  \quad \textrm{and} \quad
  M(w)_{ij}
  :=
  \int_\Omega \basis{w}_i \basis{w}_j \di{x}
  .
\end{equation}
We emphasise that both the vector $\vec{\omega}\bc{u}$ and these
matrices are of dimension $3(q+1)M_x$, the same as our solution vector, $\vec{z}^{n+1}$ (the complete
discrete approximation).  Indeed, we can rewrite our conservation laws
in terms of dot products and inner products of these vector/matrices
and $\vec{z}^{n+1}$, as
\begin{equation} \label{eqn:lkdvclaws}
  \begin{split}
    g_1(\vec{z}^{n+1}) = \vec{\omega}(u)^T \vec{z}^{n+1} & = \int_\Omega U^0 \di{x}
    ,\\
    g_2(\vec{z}^{n+1}) = \frac12 \bs{\vec{z}^{n+1}}^T M(u) \vec{z}^{n+1} & = \frac12
    \int_\Omega \bs{U^0}^2 \di{x}
    , \\
    g_3(\vec{z}^{n+1}) = \frac12 \bs{\vec{z}^{n+1}}^T M(w) \vec{z}^{n+1}
    - \frac12 \bs{\vec{z}^{n+1}}^T M(u) \vec{z}^{n+1}
    &  =
    \frac12 \int_\Omega \bs{W^0}^2 - \bs{U^0}^2 \di{x}
    .
  \end{split}
\end{equation}
We note that the ``local'' expression of these conserved quantities is that
$g_\ell(\vec{z}^{n+1})=g_\ell(\vec{z}^n)$ for $\ell=1,2,3$; however,
we have iterated these laws back to the initial data, which shows that
there values are constants, determined only by $U^0$ and $W^0$.

\revised{To} generalise the temporal discretisation used in
\eqref{eqn:fulllkdv} to arbitrarily high temporal order \revised{while preserving quadratic invariants, we use higher-order GLRK discretisations~\cite{HairerLubichWanner:2006}.}

\begin{definition}[GLRK finite element discretisation for
  lKdV] \label{def:lkdvrk} 

  Let $s$ be the number of RK stages, and
  $U^j, V^j, W^j \in \dpoly{q}$ be given for $j=0,...,n$. The approximation
  $U^{n+1}, V^{n+1}, W^{n+1}$ is given by seeking stage values
  $K^i_{U}, K^i_{V}, K^i_{W} \in \dpoly{q}$ for $i=1,...,s,$ such that
  \begin{equation} \label{eqn:lkdvrk}
    \begin{split}
      \int_\Omega \bc{
        K^i_U + \gfunc{ V_i^n}
      } \phi
      \di{x}
      & = 0
      \qquad
      \forall \phi \in \dpoly{q}, 1 \leq i \leq s
      \\
      \int_\Omega \bc{
        V_i^n - U_i^n - \gfunc{W_i^n}
      } \psi \di{x}
      & = 0
      \qquad
      \forall \psi \in \dpoly{q}, 1 \leq i \leq s
      \\
      \int_\Omega \bc{
        W_i^n - \gfunc{U_i^n}
      } \chi \di{x}
      & = 0
      \qquad
      \forall \chi \in \dpoly{q}, 1 \leq i \leq s
      ,
    \end{split}
  \end{equation}
  where $U_i^n = U^n + \dt{n}\sum_{j=1}^s a_{ij}K_U^j$,
  $V_i^n = V^n + \dt{n}\sum_{j=1}^s a_{ij}K_V^j$, and
  $W_i^n = W^n + \dt{n}\sum_{j=1}^s a_{ij}K_W^j$.
  After solving for the stage values, the solution at $t=t_{n+1}$
  is constructed through
  \revised{
  \begin{equation} \label{eqn:lkdvrk2}
    \begin{split}
      U^{n+1} =
      U^n + \dt{n} & \sum_{i=1}^s b_i K^i_U
      ,
      \qquad
      V^{n+1} =
      V^n + \dt{n} \sum_{i=1}^s b_i K^i_V
      ,
      \\
      & W^{n+1}  =
      W^n + \dt{n} \sum_{i=1}^s b_i K^i_W
      .
    \end{split}
  \end{equation}
  }
\end{definition}
\revised{%
  We note that, if $U^n$, $V^n$, and $W^n$ satisfy
  {$\int_\Omega \bc{V^n - U^n - \gfunc{W^n}} \psi \di{x} = 0$} for all
  $\psi$ and \\ {$\int_\Omega \bc{W^n - \gfunc{U^n}} \chi \di{x}=0$} for
  all $\chi$, then the second two equations in~\eqref{eqn:lkdvrk} can
  be enforced directly on the stage values, rather than on the full
  stage approximations, and these relations will also hold at
  time-step $n+1$.  These can be enforced at time-step $0$ by
  requiring that $W^0 = \gfunc{U^0}$, where $U^0$ is the given
  initial data for $U$ (typically, $U^0 = \Pi u_0$), and that
  $V^0 = U^0 + \gfunc{W^0}$.
}

When $s=1$, Definition~\ref{def:lkdvrk} is equivalent to the \CN
scheme; however, we note that the discrete solution process is not
identical, since the \CN scheme solves directly for
$U^{n+1}, V^{n+1}, W^{n+1}$, while the one-stage Gauss-Legendre (implicit midpoint)
scheme solves for the stage values, then updates $U^n, V^n, W^n$ to
get the values at the next time-step.

The conserved quantities for the higher-order scheme are unchanged and
given by \eqref{eqn:lkdv:laws}.  These constraints are enforced in
much the same way, except that the solution to the linear
system $\revised{\A} \vec{x} = \revised{\vec{f}}$ is no longer $\vec{z}^{n+1}$, the basis
coefficients of $\vec{Z}^{n+1}$.  Instead, we must accumulate the
stage values encoded in $\vec{x}$ to get $\vec{z}^{n+1}$, given by
\begin{equation} \label{eqn:lkdvreconstruction}
  \vec{z}^{n+1} = \vec{z}^n + \dt{n} \sum_{i=1}^s b_i \vec{x} \big|^i
  ,
\end{equation}
where $\vec{z}^n$ is the solution at time $t=t_n$ and
$\vec{x} \large|^i$ is the restriction of $\vec{x}$ to the $i$-th
stage. With \eqref{eqn:lkdvreconstruction} in mind, an $s$-stage GLRK
method exactly preserves the algebraic form of the constraints
\eqref{eqn:lkdvclaws}.

\subsection{Shallow water equations} \label{sec:background:swe}

We next consider a two-dimensional domain, with
$\vec{x} \in \Omega = [0,X] \times [0,Y]$ and $t\in[0,T]$. On that
domain, we consider $\vec{u}(t,\vec{x}) \in \real^2$ representing a
horizontal velocity and $\rho(t,\vec{x}) \in \real$ a proxy for
pressure, and define the linear rotating shallow water equations,
\begin{equation} \label{eq:swe1}
  \begin{split}
    \vec{u}_t + f \vec{u}^\perp + c^2 \nabla\rho & = 0 \\
    \rho_t + \nabla \cdot \vec{u} & = 0
    ,
  \end{split}
\end{equation}
where $f$ is the Coriolis parameter, $c^2$ is the gravitational
acceleration times the mean layer thickness, and $\vec{u}^\perp =
\vec{e}_3 \times \vec{u}$.  We impose doubly periodic boundary
conditions for both $\vec{u}$ and $\rho$ on $\partial\Omega$.

We use the discretisation from
\cite{CotterShipton:2012}, which preserves many of the geometric
properties of the shallow water equations on triangulations of $\Omega$.  Letting
$H(\operatorname{div},\Omega)$ be the space of vector fields,
$\vec{u}\in \left(L^2(\Omega)\right)^2$, such that $\nabla\cdot\vec{u}
\in L^2(\Omega)$, we consider approximating $\vec{u}$ in the
Raviart-Thomas space of degree $q$, denoted $\RT_q$.  (We note that,
in our notation, the lowest-order Raviart-Thomas space is $\RT_1$, and
that the basis functions for $\RT_q$ are subsets of polynomials of
degree at most $q$.)  To represent $\rho$, we use discontinuous
Lagrange elements of degree $q-1$, denoted again by $\dpoly{q-1}$
(although we note that this is now the two-dimensional discontinuous
Lagrange space on triangles).
With these, and an integration-by-parts on $\nabla\rho$,  the
semi-discretised finite-element method seeks $\vec{U} \in \RT_q$ and
$\Rho \in \dpoly{q-1}$ such that
\begin{equation} \label{eqn:swe:semi}
  \begin{split}
  \int_\Omega \vec{U}_t \cdot \phi + f \vec{U}^\perp \cdot \phi
  - c^2 \Rho \nabla \cdot \phi \di{x} & = 0
  \qquad \forall \phi \in \RT_q \\
  \int_\Omega \Rho_t \psi + \nabla \cdot \vec{U} \psi \di{x} & =
  0
  \qquad \forall \psi \in \dpoly{q-1}
  .
  \end{split}
\end{equation}
This scheme takes advantage of the discrete de Rham complex of
finite-element exterior calculus~\cite{Arnold2006, Cotter2014}, with
the property that for
$\vec{U} \in \RT_q$, we have $\nabla \cdot \vec{U} \in \dpoly{q-1}$.  We
note that~\cite{CotterShipton:2012} advocates for use of higher-order
analogues of the $\mathrm{BDFM}_1$ space in place of $\RT_q$.  This is
necessary to preserve additional important invariants in the context
of atmospheric dynamics, but not for the simpler invariants that we
consider here.

We utilise \CN for our temporal discretisation. Due to the linear
nature of the problem, the invariants will be quadratic, thus the
temporal discretisation shall be conservative.

\begin{definition}[Finite element discretization for SWE] \label{def:swe:scheme}

  Let $\vec{U}^j \in \RT_q$, $\Rho^j \in \dpoly{q-1}$ be given for
  $j=0,...,n$, then we seek $\vec{U}^{n+1}\in \RT_q$ and
  $\Rho^{n+1} \in \dpoly{q-1}$ such that
  \begin{equation} \label{eqn:swe:scheme}
    \begin{split}
      \int_\Omega \bc{ \frac{\vec{U}^{n+1}-\vec{U}^n}{\dt{n}}} \cdot
      \phi
      + f \bs{\vec{U}^{n+\frac12}}^\perp \cdot \phi
      - c^2 \Rho^{n+\frac12} \nabla \cdot \phi \di{x} & = 0
      \qquad \forall \phi \in \RT_q \\
      \int_\Omega \bc{\frac{\Rho^{n+1} - \Rho^n}{\dt{n}}} \psi + \nabla \cdot \vec{U}^{n+\frac12} \psi \di{x} & =
      0
      \qquad \forall \psi \in \dpoly{q-1}
      ,
    \end{split}
  \end{equation}
where $\vec{U}^{n+\frac12} = \frac12\left(\vec{U}^n +
\vec{U}^{n+1}\right)$ and $\Rho^{n+\frac12} = \frac12\left(\Rho^n + \Rho^{n+1}\right)$.
\end{definition}
  
This finite-element discretisation conserves a plethora of geometric
properties (cf.,~\cite{CotterShipton:2012}).  Here, we focus on the conservation laws for mass and
energy, namely that
\begin{equation}
  \int_\Omega \Rho^{n+1} \di{x} = \int_\Omega \Rho^n \di{x}
\end{equation}
and
\begin{equation}
  \frac12 \int_\Omega \abs{\vec{U}^{n+1}}^2 + c^2 \bs{\Rho^{n+1}}^2
  \di{x}
  =
  \frac12 \int_\Omega \abs{\vec{U}^n}^2 + c^2 \bs{\Rho^n}^2 \di{x}
  ,
\end{equation}
respectively. As before, we define
$\vec{Z}^{n+1} = [\vec{U}^{n+1};\Rho^{n+1}]$ and denote the solution
vector in terms of the basis coefficients as $\vec{z}^{n+1}$. As we
directly solve for $\vec{z}^{n+1}$ in the \CN discretisation, we have
$\vec{z}^{n+1}=\vec{x}$, the solution of $\revised{\A}\vec{x} = \revised{\vec{f}}$.  To reformulate the constraints, we define
the weight vector and mass matrices componentwise as
\revised{
\begin{equation}
  \omega(\rho)_i
  :=
  \int_\Omega \basis{\rho}_i \di{x}
  ,
  \qquad
  M(\vec{u})_{ij}
  :=
  \int_\Omega \basis{\vec{u}}_i \cdot \basis{\vec{u}}_j \di{x}
  ,
  \qquad
  M(\rho)_{ij}
  :=
  \int_\Omega \basis{\rho}_i \basis{\rho}_j \di{x}
  ,
\end{equation}
}%
where the linear operators are formed through decomposing the trial
and test functions in terms of their basis functions. We note that
while these operators are formulated similarly to their analogues for the linear KdV
equation, these are, in practice, quite different due to the differing
natures of the basis functions. With these operators, we
write the constraints on the underlying linear system as
\begin{subequations} \label{eqn:swe:constraints}
\begin{align}
    g_1(\vec{z}^{n+1}) & = \vec{\omega}(\rho)^T \vec{z}^{n+1} =
    \int_\Omega \Rho^0 \di{x} \\
    g_2(\vec{z}^{n+1}) & = \frac12 \left(\vec{z}^{n+1}\right)^T
    M(\vec{u}) \vec{z}^{n+1}
    + \frac{c^2}{2} \left(\vec{z}^{n+1}\right)^T M(\rho) \vec{z}^{n+1} =
    \frac12 \int_\Omega \abs{\vec{U}^0}^2 + c^2 \bs{\Rho^0}^2
    \di{x}\notag
    ,
\end{align}
\end{subequations}
where the right-hand sides are, again, constants determined by the
initial conditions.

\subsection{Heat equation} \label{sec:background:heat}

Let $u(t,\vec{x})\in\real$ for
$\vec{x}\in \Omega = [0,1] \times [0,1]$ and $t\in[0,1]$, then the
heat equation may be written as
\begin{equation}
  \label{eqn:heat}
  u_t - \Delta u = 0
  ,
\end{equation}
where, for simplicity, we employ Neumann boundary conditions in
space. We consider a standard discretisation
of the heat equation using conforming Lagrangian finite elements, writing
the space of $q$-th order triangular Lagrangian elements over
the domain $\Omega$ as $\bpoly{q}$. The spatial finite-element
semi-discretisation of~\eqref{eqn:heat}
is then given by seeking $u \in \bpoly{q}$ such that
\begin{equation}
  \int_\Omega U_t \phi +  \nabla U \cdot \nabla \phi \di{x} = 0
  \qquad
  \forall \phi \in \bpoly{q}
  ,
\end{equation}
subject to some initial conditions $U(0,\vec{x}) = \Pi
u_0(\vec{x})$. We again discretise
temporally with \CN leading to the following fully discrete method.

\begin{definition}[Finite element discretization for the heat
  equation]

  Let $U^j \in \bpoly{q}$ be given for $j=0,...,n$, then we seek
  $U^{n+1}$ such that
  \begin{equation} \label{eqn:heat:scheme}
    \int_\Omega \bc{\frac{U^{n+1}-U^n}{\dt{n}}} \phi
    + \nabla U^{n+\frac12} \cdot \nabla \phi \di{x} = 0 
    \qquad
    \forall \phi \in \bpoly{q}
    ,
  \end{equation}
  where $U^0 = \Pi u_0(\vec{x})$ and $U^{n+\frac12} = \frac12\left(U^n
  + U^{n+1}\right)$.
\end{definition}

The heat equation respects conservation of mass and a dissipation law,
and its finite-element approximation preserves discrete counterparts
of these; namely, the mass conservation law
\begin{equation}
  \int_\Omega U^{n+1} \di{x}
  =
  \int_\Omega U^n \di{x}
\end{equation}
and the dissipation law
\begin{equation}
\frac{1}{2\dt{n}}\int_\Omega \bs{U^{n+1}}^2 - \bs{U^{n}}^2 \di{x}=
-\int_\Omega \nabla U^{n+\frac12}\cdot \nabla U^{n+\frac12} \di{x},
\end{equation}
which can be rewritten as
\begin{equation}
  \frac12 \int_\Omega \bs{U^{n+1}}^2
  + \frac{\dt{n}}{4} \bs{ \nabla U^{n+1}}^2
  + \frac{\dt{n}}{2} \nabla U^{n+1} \cdot \nabla U^n
  \di{x}
  =
  \frac12 \int_\Omega \bs{U^n}^2
  - \frac{\dt{n}}{4} \bs{\nabla U^n}^2
  \di{x}
  .
\end{equation}

To reformulate mass conservation and energy dissipation in terms of
the underlying linear system, we require the weight vector and mass
matrix
\begin{equation}
  \omega_i = \int_\Omega \basis{}_i \di{x}
  \quad
  \textrm{and}
  \qquad
  M_{ij} = \int_\Omega \basis{}_i \basis{}_j \di{x}
  ,
\end{equation}
in addition to the stiffness matrix $L$ defined such that
\begin{equation}
  L_{ij}
  =
  \int_\Omega
  \nabla \basis{}_i \cdot \nabla \basis{}_j
  \di{x}
  .
\end{equation}
We note that, while the mass conservation law can again be formulated
in terms of only the initial condition, the dissipation law is more
naturally expressed in terms of time-steps $n$ and $n+1$.  From the
discrete solution $\vec{z}^{n+1}$ to the \CN discretisation, we can
express the constraints as
\begin{subequations} \label{eqn:heat:constraints}
\begin{align}
    g_1(\vec{z}^{n+1})
    & =
    \vec{\omega}^T \vec{z}^{n+1}  = \int_\Omega U^0 \di{x} \\
    g_2(\vec{z}^{n+1}) &= \frac12 \left(\vec{z}^{n+1}\right)^T M \vec{z}^{n+1}
    + \frac{\dt{n}}{4} \left(\vec{z}^{n+1}\right)^T L \vec{z}^{n+1}
    + \frac{\dt{n}}{2} \left(\vec{z}^{n+1}\right)^T L \vec{z}^n
     =
    \frac12 \left(\vec{z}^n\right)^T M \vec{z}^n
     - \frac{\dt{n}}{4} \left(\vec{z}^n\right)^T L \vec{z}^n
    .\notag
  \end{align}
\end{subequations}
We note that $g_2(\vec{z}^{n+1})$ is no longer a constraint in the
form of $g_2(\vec{z}^{n+1}) = g_2(\vec{z}^{n})$, since this represents
a dissipation law, not a conservation law.  While we write the
algorithms below as if the constraints are always conservation laws,
we note that the algorithm can, in fact, apply to any relation
$G(\vec{z}^{n+1},\vec{z}^n) = 0$.  To express this concisely in what
follows, we consider a generic linear algebra problem,
$\revised{\A}\vec{x}=\revised{\vec{f}}$, coming from one of these discretisations, along
with a list of $c$ constraints to be satisfied by the discrete solution,
$\vec{x}$, of the form $g_i(\vec{x}) = v_i$ for $1\leq i \leq c$,
noting that all of the constraints above can be written in this form,
for both the \CN discretisations (where $\vec{x} = \vec{z}^{n+1}$) and
the higher-order GLRK discretisations (where $\vec{z}^{n+1}$ is
assembled from $\vec{x}$, as in~\eqref{eqn:lkdvreconstruction}).

\section{Structure-preserving Krylov methods}\label{sec:methodology}

Krylov methods~\cite{AGreenbaum_1997, YSaad_2003a} are iterative
algorithms for approximating solutions to linear systems
$\revised{\A}\vec{x}=\revised{\b}$, by generating approximations to $\vec{x}$ in the
Krylov space given by
\begin{equation}
\mathcal{K}_{\ell}(\revised{\A},\vec{r}_0) = \text{span}\left\{\vec{r}_0,
  \revised{\A}\vec{r}_0, \ldots, \revised{\A}^{\ell-1}\vec{r}_0\right\},
\end{equation}
where $\vec{r}_0 = \revised{\b}-\revised{\A}\vec{x}_0$ is the initial residual
associated with an initial guess, $\vec{x}_0$, for $\vec{x}$. There
are multiple (overlapping) classifications of Krylov methods. Some
aim to minimise some norm of the error associated with $\vec{x}_\ell
\in \mathcal{K}_{\ell}(\revised{\A},\vec{r}_0)$, while others aim to enforce
some orthogonality relationship between that error and the Krylov
space itself.  Most commonly used techniques can be considered in the
first class.  Amongst these, we can further distinguish between those
that make use of short-term recurrences (e.g., conjugate gradients for
symmetric and positive-definite matrices, or MINRES for symmetric
systems), and those that form and store an explicit basis (usually
orthogonal in some inner product) for the Krylov space, such as GMRES.
While it may be possible to adapt our approach to any method that
stores an explicit basis, we focus on adapting flexible GMRES (FGMRES)
to preserve conservation laws of the approximate solutions generated
by an iterative method.

We note that there exist two commonly used variants of the GMRES
methodology, that differ in how preconditioners are incorporated.
Standard GMRES aims to choose $\vec{x}_\ell$ to minimise the residual,
$\revised{\b}-\revised{\A}\vec{x}_\ell$, over the affine space $\vec{x}_\ell =
\vec{x}_0 + \delta\vec{x}_\ell$, where $\delta\vec{x}_\ell \in
\mathcal{K}_{\ell}(\revised{\A},\vec{r}_0)$.  The standard approach is to
use the Arnoldi algorithm to form an orthogonal basis for
$\mathcal{K}_{\ell}(\revised{\A},\vec{r}_0)$, expressed as the range of the $sd \times \ell$
matrix, $Q_\ell$ (recall that $\revised{\A} \in \mathbb{R}^{sd\times
  sd}$, where $s$ is the number of Runge-Kutta stages, and $d$ is the
dimension of the spatial discretisation).  Then, defining the
$(\ell+1)\times \ell$
Hessenberg matrix, $H_\ell$, associated with the Arnoldi basis in
$Q_\ell$, whose first column is $\vec{r}_0$, the problem of minimising
$\|\revised{\b} - \revised{\A}\vec{x}_\ell\|$ for $\vec{x}_\ell = \vec{x}_0 + Q_\ell
\vec{y}_\ell$ is equivalent to that of minimising
$\|\beta\vec{e}_{\ell+1} - H_\ell \vec{y}_\ell\|$, for $\beta =
\|\vec{r}_0\|$ and $\vec{e}_{\ell+1}$ the unit vector of length $\ell+1$
with first entry 1 and all other entries zero.  Standard GMRES uses
Householder transformations to iteratively compute a representation of the QR
factorisation of $H_\ell$ (updating that of $H_{\ell-1}$) that allows
us to compute the value of the norm at each step, then solve a
triangular system for $\vec{y}_\ell$ once a stopping criterion is
satisfied (see~\cite{YSaad_2003a} for full details).  Left or right
preconditioners can be incorporated by replacing the matrix $\revised{\A}$ with
$\precon{}\revised{\A}$ or $\revised{\A}\precon{}$, respectively, in the Arnoldi
algorithm, to find $\vec{y}_k$ that minimises either
\begin{equation}
  \left\|\precon{}\revised{\b} - \precon{}\revised{\A}\left(\vec{x}_0 + Q_\ell
  \vec{y}_\ell\right)\right\|
\end{equation}
for left preconditioning or
\begin{equation}
  \left\|\revised{\b} - \revised{\A}\left(\vec{x}_0 + \precon{}Q_\ell
  \vec{y}_\ell\right)\right\|
\end{equation}
for right preconditioning.

For many problems, right preconditioning is preferable to left
preconditioning, as it does not change the norm of the underlying
minimisation problem, just the Krylov space over which we minimise,
becoming $\mathcal{K}(\revised{\A}\precon{},\vec{r}_0)$.  When
right-preconditioning with classical GMRES, however, we require
$\ell+1$ applications of $\precon{}$ to compute $\vec{x}_\ell$,
with $\ell$ used to compute the basis for the Arnoldi space, and one
more to compute $\vec{x}_\ell = \vec{x}_0 + \precon{}Q_\ell
\vec{y}_\ell$ once $\vec{y}_\ell$ is known.  The flexible GMRES
(FGMRES) variant removes this ``extra'' application of the
preconditioner in favour of doubling the vector storage, to store both
the standard orthogonal basis for the Krylov space, $Q_\ell$, and its
preconditioned form, $Z_\ell = \precon{}Q_\ell$ (computed and stored
columnwise within the preconditioned Arnoldi iteration).  Then, we can
compute $\vec{x}_\ell = \vec{x}_0 + Z_\ell \vec{y}_\ell$ once
$\vec{y}_\ell$ is known.  Unless the target computation is severely
memory-limited, using FGMRES is generally preferred, particularly when
the cost of preconditioner application is high (as it is in many cases
when $\revised{\A}$ comes from the discretisation of coupled systems of PDEs).
This also offers the flexibility to use different preconditioners at
different inner iterations, with $\precon_\ell$ applied at step
$\ell$, although we do not exploit this feature in what follows. For
completeness, we sketch the FGMRES algorithm in
Algorithm \ref{alg:gmres}, albeit omitting the important details of the QR
factorisation of $H_\ell$.
\begin{algorithm2e}[h]
  \DontPrintSemicolon
  
  \KwInput{$\revised{\A}$, $\revised{\b}$, $\vec{x}_0$\tcp*{Matrix, right-hand side,
      initial guess}}
  \nonl \myinput{$\ell_{\text{max}}$ \tcp*{Maximum number of iterations}}
  \nonl \myinput{$\precon$ \tcp*{Preconditioner (optional, defaults to identity matrix), may depend on iteration, $\ell$}}
  \nonl \myinput{$\tol$ \tcp*{Convergence tolerance}}
  
  \KwOutput{$\vec{x}_\ell$ \tcp*{Returns an
      approximation of the linear system}}
  
  $\vec{r}_0 = \revised{\b} - \revised{\A}\vec{x}_0$

  $\beta = \| \vec{r}_0 \|$
  
  $\vec{q}_1 = \vec{r}_0 / \beta$

  \For(\tcp*[f]{Loop over iterations}){$\ell=1,...,\ell_{\text{max}}$}
  {
    $\vec{z}_\ell = \precon_\ell \vec{q}_\ell$
    
    $\vec{q} = \revised{\A} \vec{z}_\ell$

    \For(\tcp*[f]{Arnoldi iteration}){$i=0,...,\ell+1$}
    {
      $h_{i\ell} = \vec{q}_i^T \vec{q}$ \;
      $\vec{q} = \vec{q} - h_{i\ell} \vec{q}_i$
    }
    $h_{\ell+1,\ell} = \| \vec{q} \|$ \;
    \If{$h_{\ell+1,\ell} \neq 0$}
    {
      $\vec{q}_{\ell+1} = \vec{q} / h_{\ell+1,\ell}$
    }
    
    $\vec{y}_\ell = \underset{\vec{y}}{\min} \| \beta e_{\ell+1} - H_\ell\vec{y} \|$ \tcp*{Find
      minimiser in residual over Krylov
      space} \label{alg:gmres:solve}
    \If{$\| \beta e_{\ell+1} - H_\ell\vec{y}_\ell \| < \tol$}
       {
         Break \label{alg:gmres:check}
       }
  }
  $\vec{x}_\ell = \vec{x}_0 + Z_\ell \vec{y}_\ell$ \tcp*{Assemble solution at
      current step}\label{alg:gmres:assemble_x}
   
\caption{FGMRES \label{alg:gmres}}
\end{algorithm2e}

\revised{
\begin{remark}[Restarting]\label{rem:restarting}
  For simplicity, we do not include
  restarting in any of the algorithms presented here, as it does not play a role in any of the numerical results.  In many settings, however, particularly when the choice of optimal preconditioners is difficult, restarting is an important factor in limiting the memory consumption of GMRES-type algorithms.  In Algorithm~\ref{alg:gmres}, restarting would take the form of an outer loop over all of the presented algorithm, looping until either convergence is achieved or some maximum number of restarts is reached.  If convergence is not reached before $\ell_{\text{max}}$ steps, we compute $\vec{x}_{\ell_{\text{max}}}$ in Step~\ref{alg:gmres:assemble_x} and use it to overwrite $\vec{x}_0$ to restart the outer loop.
\end{remark}
}%

\subsection{A prototype for constrained GMRES} \label{sec:pcgmres}

The core idea of this paper is to constrain FGMRES such that desired
properties of the discretisation are satisfied. This is achieved
through modifying Line \ref{alg:gmres:solve} in Algorithm
\ref{alg:gmres}. Instead of directly minimising the residual over the
Krylov subspace, we pose a constrained minimisation problem,
ensuring that the reconstructed solution,
$\vec{x}_\ell = \vec{x}_0 + Z_\ell \vec{y}_\ell$, associated with the
(constrained) minimiser, $\vec{y}_\ell$, satisfies the desired set of
constraints.
\revised{That is, we replace Line~\ref{alg:gmres:solve} with a minimisation problem of the form}
\begin{equation}\label{eq:cgmres:solve}
  \begin{split}
    \vec{y}_\ell & = \underset{\vec{y}\in \mathcal{Y}_{\ell}}{\min} \left\| \beta e_{\ell+1} -
      H_\ell\vec{y} \right\| \\
    & \qquad\qquad \text{ for } \mathcal{Y}_{\ell} =
    \left\{ \vec{y}\in\mathbb{R}^\ell \middle| g_i\left(\vec{x}_0 + Z_\ell \vec{y}\right)
      = v_i \text{ for } 0\leq i \leq \min(\ell-1,c)\right\},
  \end{split}
\end{equation}
\revised{where we consider an additional input to the algorithm of $\clist = \{g_1, g_2, \ldots, g_c\}$, a list of constraints, each of form $g_\ell(\vec{x}) = v_\ell$.}

In the general case (when $\ell > c$), this can be thought of in terms
of a nonlinear optimisation problem in $\ell+c$ variables, with the
$c$ constraints imposed via Lagrange multipliers.  In practice, we use
the Python library \verb|SciPy|, using the \verb|optimize.minimize|
function.  This employs a \revised{sequential quadratic programming algorithm,} \verb|SQSLP|\revised{,
which wraps the original (Fortran) implementation by Kraft~\cite{DKraft_1988}.  In preliminary work, we also considered a Byrd-Omojokun Trust-Region SQP method~\cite[Section 15.4.2]{Conn:2000}, as implemented in the} \verb|trust-constr| \revised{algorithm, but found this to be significantly less efficient when timing the results, likely because it does not allow us to directly provide gradients of the constraints}.  We note that this is
substantially less efficient than the updating QR factorisation
typically used for the Hessenberg matrix, $H_\ell$, in standard
(F)GMRES.  While it may be possible to optimize this solve in the case
of only linear constraints (where the resulting first-order \revised{optimality} condition
is an $(\ell+c)\times(\ell+c)$ linear system), we are unaware of any
optimisations possible in the general case.

A crucial limitation of this methodology is that the vector subspace
for the optimisation in \revised{\eqref{eq:cgmres:solve}} must be rich
enough for a minimiser subject to the desired constraints to exist,
which is not guaranteed for arbitrary choices of $\vec{x}_0$,
particularly when $\ell$ is small.  \revised{Equation~\eqref{eq:cgmres:solve}}
mitigates this possibility by limiting the number of constraints
enforced at step $\ell$ to be no more than $\ell$.  Of course, there
is no guarantee that any given Krylov subspace will contain solutions
which satisfy any given constraint; however, the constraints that we
seek to impose here generally arise ``naturally'' within the system we
are solving (in the sense that they are guaranteed to be satisfied by
the exact solution of the underlying linear system), which gives some
confidence that a constrained minimiser should exist, at least for
sufficiently large $\ell$. With this in mind, we increase the number
of constraints enforced as we loop through the iterations, starting
with zero constraints and adding a new constraint each iteration. This
procedure could be modified to allow for more or fewer constraints to
be enforced in the early iterations, but preliminary experiments
showed that adding one per iteration was sufficient for the majority
of the test cases considered.  We note that \revised{directly substituting \eqref{eq:cgmres:solve} for Line~\ref{alg:gmres:solve} in Algorithm~\ref{alg:gmres}}
implicitly assumes that the stopping criterion will not be satisfied until
$\ell \geq c$.  While this was true in practice in all experiments
considered herein, it could easily be made an explicit condition on
Line~\ref{alg:gmres:check} if needed.

\subsection{Practical implementation of CGMRES} \label{sec:alt}

Even when the Krylov space is rich enough such that a constrained
minimiser exists \revised{for \eqref{eq:cgmres:solve}}, the computational cost of finding such
minimiser is substantially higher than that of the usual updating QR
factorisation used in (F)GMRES.  Indeed, performing the constrained
solves can be particularly expensive when the Krylov space is small,
as significant effort can be required to find (or fail to find) a constraint-satisfying
solution.  Taken together, these suggest that a practical
implementation of the constrained GMRES approach should be one where
the use of constrained solves is avoided unless we are close enough to
solution that the effort is worthwhile.  One excellent indicator of
whether or not this computational cost is worth expending is the size
of the \textit{unconstrained} residual, noting that we still use a
residual-based stopping tolerance \revised{on constrained GMRES} and
that the residual in the constrained solution cannot be smaller than
that in the unconstrained solution.  Since the norm of the
unconstrained residual can be computed at relatively low cost, this
leads to an effective strategy for avoiding the nonlinear solve until
we are close to convergence.  \revised{In fact, a regular FGMRES routine could be used until this point, giving residual norm estimates from the standard approach using QR factorisation of the Hessenberg systems at no additional cost.}

This approach is presented in Algorithm~\ref{alg:cgmres2}.  Here, in
addition to the standard residual-norm stopping parameter, $\tol$, we
require a second tolerance, $\Epsilon$, to determine whether we are
``close enough'' to convergence to make enforcing the constraints in
the solve worthwhile.  An important practical question is, of course,
how to choose the second parameter, $\Epsilon$.  We note that we
choose to structure the algorithm so that the decision to impose
constraints is based on already computed information, at step
$\ell-1$.  Thus, we view $\Epsilon$ as being how close to convergence
we expect to be at the iteration \textit{prior} to satisfying the
convergence criterion, so that we aim to perform only a single
constrained solve.  Thus, $\Epsilon/\tol$ should be approximately the
residual reduction per iteration expected from the unconstrained (but
preconditioned) iteration, so that we only execute the solve on
Line~\ref{alg:cgmres2:solve} once (assuming we can find a
constraint-satisfying solution with similar overall residual
reduction).  For the systems considered here, we generally take
$\Epsilon = 10\tol$, although we note that different choices could
easily be justified, and that the interplay between constraints and
residual minimisation in any particular problem may drive other
choices (in particular in cases where the constrained residual norm is
much larger than its unconstrained counterpart).  
  
\begin{algorithm2e}[h]
  \DontPrintSemicolon
  
  \KwInput{$\revised{\A}$, $\revised{\b}$, $\vec{x}_0$\tcp*{Matrix, right-hand side,
      initial guess}}
  \nonl \myinput{$\ell_{\text{max}}$ \tcp*{Maximum number of iterations}}
  \nonl \myinput{$\precon$ \tcp*{Preconditioner (optional, defaults to
      identity matrix), may depend on iteration, $\ell$}}
  \nonl \myinput{$\tol$ \tcp*{Convergence tolerance}}
  \nonl \myinput{$\clist = \{g_1, g_2, \ldots, g_c\}$  \tcp*{List of
      constraints of form $g_\ell(\vec{x}) = v_\ell$}}
  \nonl \myinput{$\Epsilon$ \tcp*{Tolerance to check constraint
      enforcement, defaults to $10\tol$, must have $\tol \leq \Epsilon$}}
  
  \KwOutput{$\vec{x}_\ell$ \tcp*{Returns an
      approximation of the linear system}}
  
  $\vec{r}_0 = \revised{\b} - \revised{\A} \vec{x}_0$

  $\beta = \| \vec{r}_0 \|$
  
  $\vec{q}_1 = \vec{r}_0 / \beta$

  \For(\tcp*[f]{Loop over iterations}){$\ell=1,...,\ell_{\text{max}}$}
  {
    $\vec{z}_\ell = \precon_\ell \vec{q}_\ell$
    
    $\vec{q} = \revised{\A}\vec{z}_\ell$

    \For(\tcp*[f]{Arnoldi iteration}){$i=0,...,\ell+1$}
    {
      $h_{i\ell} = \vec{q}_i^T \vec{q}$ \;
      $\vec{q} = \vec{q} - h_{i\ell} \vec{q}_i$
    }
    $h_{\ell+1,\ell} = \| \vec{q} \|$ \;
    \If{$h_{\ell+1,\ell} \neq 0$}
    {
      $\vec{q}_{\ell+1} = \vec{q} / h_{\ell+1,\ell}$
    }
    
    \uIf{$\| \beta e_{\ell} - H_{\ell-1}\vec{y}_{\ell-1} \|> \Epsilon$ \KwAnd $\ell<\ell_{\text{max}}$}
    {
      $\vec{y}_\ell = \underset{\vec{y}}{\min} \| \beta e_{\ell+1} - H_\ell\vec{y}\|$
      \tcp*{If far from tolerance and $\ell\neq \ell_{\text{max}}$,
        use FGMRES} \label{alg:cgmres2:redirect}
    }
    \Else{
    $\vec{y}_\ell = \underset{\vec{y}}{\min} \| \beta e_{\ell+1} -
    H_\ell\vec{y} \|$ subject to
    the constraints $\{g_1,\ldots,g_{c}\}$ \tcp*{If close
      to tolerance or $\ell = \ell_{\text{max}}$, use constrained
      solve} \label{alg:cgmres2:solve}
    \revised{
      \If{\textrm{constrained solve fails}}
      {
        Go to line \ref{alg:cgmres2:redirect}}\label{alg:cgmres2:check}
    }
    \If{$\| \beta e_{\ell+1} - H_\ell\vec{y}_\ell \| < \tol$}
    {
      Break
    }
    }
  }
  $\vec{x}_\ell = \vec{x}_0 + Z_\ell \vec{y}_\ell$ \tcp*{Assemble solution at
      current step}

\caption{Optimised constrained GMRES \label{alg:cgmres2}}
\end{algorithm2e}

\revised{%
  \begin{remark}[Restarting with constraints] \label{rem:restarting_cgmres}
As above, we do not explicitly include restarting in Algorithm~\ref{alg:cgmres2}, but note that it could be included as an additional ``outer'' loop.  However, we also note a potential complication when restarting.  As noted above, it can be difficult to enforce constraints when the Krylov space is too small.  Thus, when we are close to convergence, it may be quite disadvantageous to restart, since we throw away a ``rich'' Krylov space over which we could effectively enforce the constraints and restart from a low-dimensional Krylov space, where enforcing constraints may be quite difficult.  This issue did not impact the numerical results reported below in~\cref{sec:examples}; however, we expect that further attention to these details may be needed if the CGMRES methodology is applied in combination with problems and preconditioners for which restarting is commonplace.  As written, Algorithm~\ref{alg:cgmres2} suggests imposing the constraints at the end of each restart cycle (the outer loop that is implicitly there) and again at convergence in the inner loop.  In practice, it is probably most useful to only enforce constraints once one is near convergence, provided the Krylov space is rich enough to expect success.
  \end{remark}
}

\revised{
  \begin{remark}[Additional computational cost] \label{rem:cost}

    While the exact cost of including and enforcing the constraints is highly problem dependent, we
    generally expect these to increase the cost of the algorithm. Within this
    work, we exclude problem-specific optimisations and focus on
    constraints that can be expressed in the form
    \begin{equation} \label{eqn:constraints}
      \vec{x}_\ell^T \mathcal{M} \vec{x}_\ell + \vec{x}_\ell^T \vec{v} + c
      = 0
      ,
    \end{equation}
    for $\mathcal{M}\in\mathbb{R}^{sd \times sd}$, $\vec{v}\in\mathbb{R}^{sd}$ and
    constant $c$, where $sd$ is the dimension of the underlying linear system $\mathcal{A} \vec{x} = \vec{f}$. We note here that
    \eqref{eqn:constraints} is the general form of a conserved
    quantity for linear problems, but that more general constraints might be of interest for broader classes of PDEs. As written in~\cref{eqn:constraints}, the evaluation of this constraint has a cost of
    $\mathcal{O}(sd)$, assuming a sparse constraint matrix, $\mathcal{M}$.  However, since $\vec{x}_\ell = \vec{x}_0 + Q_\ell\vec{y}_{\ell}$ for $\vec{y}_\ell \in \mathbb{R}^\ell$, with $\ell \ll sd$, we can greatly reduce this cost, by precomputing various products with $\vec{x}_0$ and $Q_\ell$.  Rewriting~\eqref{eqn:constraints}, we have
    \begin{equation}\label{eqn:reduced_constraints}
      \vec{y}_\ell^T\left(Q_\ell^T\mathcal{M}Q_\ell\right)\vec{y}_\ell + \vec{y}_\ell^T\left(2Q_\ell^TM\vec{x}_0 + Q_\ell^T\vec{v}\right) + \left(\vec{x}_0^T\mathcal{M}\vec{x}_0 + \vec{x}_0^T\vec{v} + c \right) = 0,
      \end{equation}
      noting that, in general, $Q_\ell^T\mathcal{M}Q_\ell$ will be a dense matrix even if $\mathcal{M}$ is sparse.  Thus, we can evaluate the reduced constraint at a cost of $\mathcal{O}(\ell^2)$, for computing the inner product in the first term, plus an $\mathcal{O}(\ell)$ cost for the dot product in the second.  The cost of evaluating the Jacobian of the constraint is the same.  However, in order to achieve this, we must precompute several terms.  Multiplying the sparse matrix $\mathcal{M}$ by the dense matrix $Q_\ell$ has $\mathcal{O}(sd\ell)$ cost ($\mathcal{O}(1)$ cost for each entry in the resulting $sd \times \ell$ matrix), so computing $Q_\ell^T\mathcal{M}Q_\ell$ has cost $\mathcal{O}(sd\ell^2)$, and this is the dominant cost in precomputing the terms in the reduced constraint in~\eqref{eqn:reduced_constraints}.  For perspective, this is the same cost as computing the Hessenberg system for the $\ell$-dimensional Krylov space in the Arnoldi algorithm.  In particular, when $\ell$ is small, this cost is typically dominated by the cost of the matrix-vector product and preconditioner applications in GMRES.
    We report details on practical (measured) computational costs in Section~\ref{sec:examples}, showing that this cost is not trivial but, at the same time, does not dominate the computation.
    
  \end{remark}
}

\begin{remark}[Nonlinear iterative solvers] \label{rem:nonlinear}

  Extending this approach to nonlinear PDEs is slightly more invasive
  but still feasible.  Here, we envision a Newton-Krylov approach,
  where the nonlinear equation for $\vec{z}^{n+1}$ is linearised via
  Newton's method, and the successive linearisations are solved using
  FGMRES.
  \revised{In this setting, the solution of the linear system $\A\vec{x} = \b$ represents the update to the current approximate solution of the nonlinear system, and a standard Newton linesearch method, used to improve robustness and efficiency of the nonlinear solve, updates the current approximation by computing $\vec{z}^{n+1} = \vec{z}^n + \omega \vec{x}$, with parameter $\omega$ chosen to minimize (in some sense) the nonlinear residual norm along the line defined by $\vec{z}^n$ and $\vec{x}$.  Since this weight is not known before we solve for $\vec{x}$, it is difficult to compute updates that} ensure that the solution of
  the nonlinear equation stays on the constraint.

  While we do not implement and explore such methods here, we note
  that a similar strategy to the above may be mimicked in this case.
  Close to (nonlinear) convergence, it is typical to see \revised{the weight chosen by} linesearch
  algorithms approach unit step size, $\omega \rightarrow 1$.  Thus, it is possible to
  construct a constrained Newton-Krylov method that iterates using
  standard Newton iterations (with inexact solves, such as using the
  Eisenstat-Walker criteria~\cite{eisenstat1996choosing}, and
  linesearches) until suitably close to convergence such that fixing the
  step size to be one is reasonable.  Then, the constrained GMRES
  algorithm can be applied to solve subsequent linearisations (ideally
  only the final linearisation) with constraints adapted to be posed
  on the (updated) nonlinear solution\revised{, now of the form $\vec{z}^{n+1} = \vec{z}^n + \vec{x}$}, rather than \revised{directly on} the solution\revised{, $\vec{x}$,} of the
  Jacobian system.
\end{remark}

\section{Numerical experiments}\label{sec:examples}

Here we apply the methodology discussed in Section~\ref{sec:methodology} for
the problems and discretisations discussed in Section~\ref{sec:background},
highlighting the potential benefits and pitfalls when incorporating
geometric constraints into linear solvers. Our implementation can be
found in \cite{self:geosolvecode}, making use of the Firedrake and
Irksome libraries for the spatial and temporal
discretisations~\cite{Firedrake:2017, Irksome:2021}. The specific
versions of the Firedrake software used are recorded
here~\cite{zenodo/Firedrake-20221208.0}. \revised{All results
  presented here, including timings, have been computed on a single
  Xeon 2.4 GHz CPU core, with access to 128 GB of memory.}

\subsection{Linear KdV}\label{sec:examples:lkdv}

Recall the discretisation of the linear KdV equation discussed
in Section~\ref{sec:background:lkdv}. We begin by solving the finite element
approximation over a single time step subject to the initial
condition
\begin{equation} \label{eqn:lkdv:initial}
  U^0 = \sin{\alpha x} + 1
  ,
\end{equation}
where $\alpha=\frac{\pi}{5}$. Note that we specifically choose $U^0$
to not have zero mean, as zero-mean solutions can achieve conservation
of mass automatically when using standard FGMRES with a zero-mean
initial guess, \revised{as shown in \cite[Section
  3.3]{BirkenLinders:2021}.} 

\begin{figure}[h!]
  \caption{A comparison between CGMRES as \revised{proposed in Section~\ref{sec:pcgmres}} and FGMRES with no preconditioner. Here, the linear system we solve
    corresponds to the finite element scheme for linear KdV
    \eqref{eqn:fulllkdv} where ${\dt{}}=0.01$, ${\M}=50$ and the spatial degree is
    $q=1$. We initialise $U^0$ as given by \eqref{eqn:lkdv:initial}
    and solve for $(U^1, V, W^1)$. We take an initial guess for both
    linear solvers of $\vec{x}_0=\vec{0}$ and solve CGMRES subject to
    the constraints \eqref{eqn:lkdvclaws}, varying the order in
    which the constraints are enforced.}
  \label{fig:lkdv}

  \includegraphics[width=0.9\textwidth]{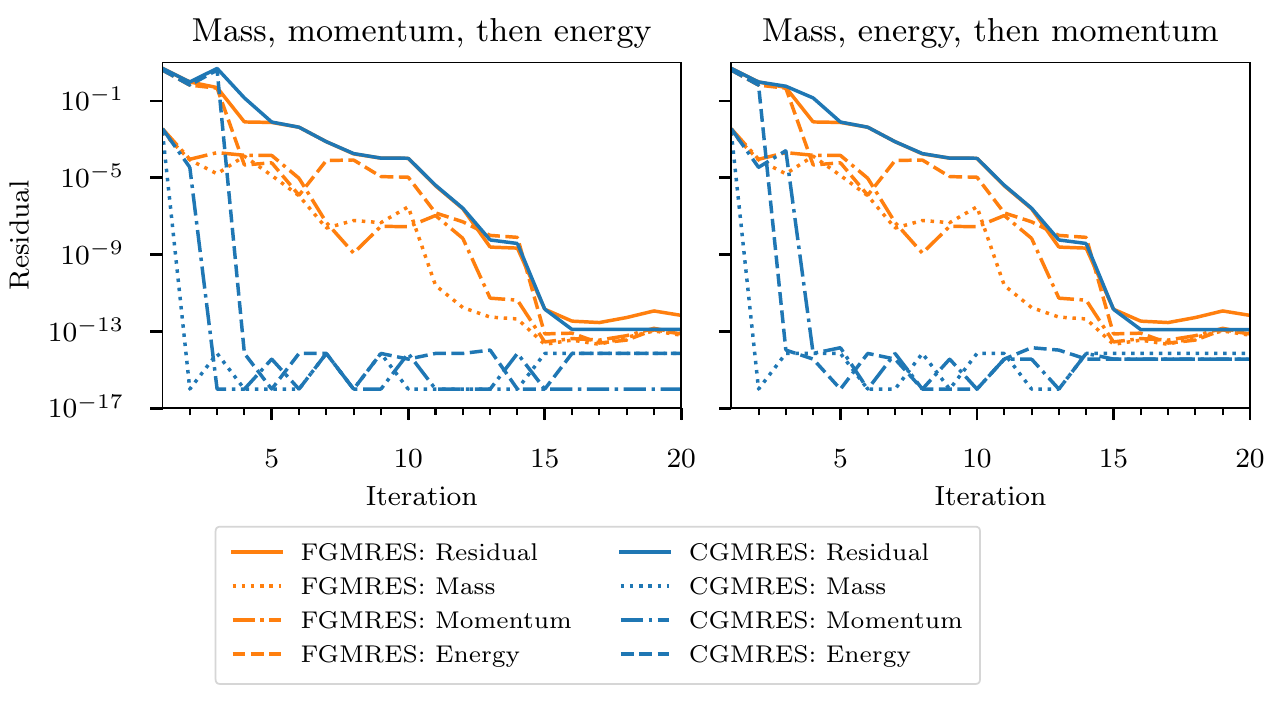}

\end{figure}

In Figure~\ref{fig:lkdv}, we compare solutions generated using the
CGMRES algorithm \revised{proposed in Section~\ref{sec:pcgmres}} against those generated with classical
FGMRES, exploring the effects of the ordering in which the constraints are enforced with the \CN
temporal discretisation \eqref{eqn:fulllkdv}.  Here, we run both
solvers for 20 iterations, to compare histories of residual norm
reduction and constraint satisfaction.  In both cases
considered, we notice that the
residual norm when using CGMRES is slightly larger than that with
FGMRES for a few early iterations, \revised{while that with FGMRES is slightly larger than that with CGMRES} near convergence.  \revised{At convergence,} we note that residuals are well below discretisation
error for this system \revised{for both methods}.  In contrast, the conserved
quantities are quickly enforced to the level of machine precision for
the CGMRES solution while
they decrease at the same rate as the residual norm for classical FGMRES. Interestingly, we note
that if we constrain mass, then energy, and then momentum, our
algorithm successfully finds constrained solutions at every
step, as shown at right of Figure~\ref{fig:lkdv}, \revised{with a steady decrease in the CGMRES residual norm and}
machine-precision errors in mass from the second iteration onwards, in
energy from the third iteration onwards, and in momentum from the
fourth iteration onwards. However, if we constrain mass, then momentum, then energy,
the constrained solve on the third iteration \revised{leads to an increase in the residual norm when we impose all three constraints for the first time}.  We expect that this is due to the Krylov space not yet being rich enough to
allow a \revised{good} solution that satisfies \revised{all three
  constraints}. As we
continue to iterate, however, and the space over which we minimise grows, 
the constraints are more easily satisfied by the approximate
solution \revised{and the residual norm for CGMRES decreases as expected}. Note that if, instead of running a fixed number of iterations, we terminate
the linear solvers at a typical residual-norm stopping tolerance of $\tol=10^{-6}$, CGMRES performs comparably to
FGMRES in both cases in terms of the residual norm, and allows us to enforce all
constraints to the level of machine precision.

\begin{figure}[h!]
  \caption{The deviation in mass, momentum and energy over time for
    the finite element scheme for linear KdV \eqref{eqn:fulllkdv} where
    ${\dt{}}=0.01$, ${\M}=50$, $T=1$ and the degree is $q=1$, using the
    linear solvers FGMRES and CGMRES as described by Algorithm
    \ref{alg:gmres} and \ref{alg:cgmres2}, respectively. Here, we fix
    $\tol=10^{-6}$ and constrain CGMRES by \eqref{eqn:lkdvclaws} with $\Epsilon=10\tol$. We
    initialise $U^0$ as given by \eqref{eqn:lkdv:initial} and solve
    iteratively for $(U^{\N}, V^{\N}, W^{\N})$. \revised{At left, we present results for all three constraints using a zero initial guess for the linear solvers at each time step.  At right, we present results for all three constraints using the solution from the previous time step as the initial guess at each time step, but only enforcing constraints on momentum and energy in CGMRES.}}
   \label{fig:lkdvevolve}
   \includegraphics[width=0.9\textwidth]{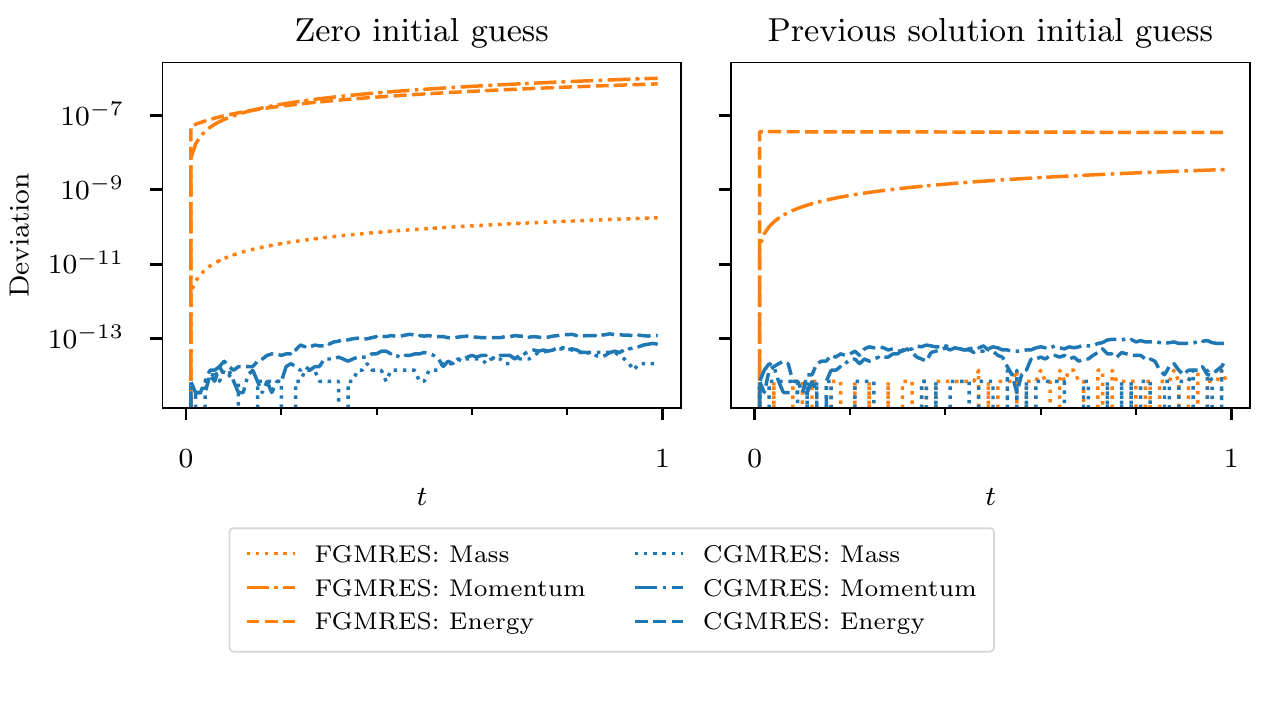}
\end{figure}

In addition to understanding the convergence of the iterative linear
solvers over a single solve, it is crucial to understand how these
algorithms impact simulations as they propagate approximate solutions over
time. With this in mind, we use a standard stopping tolerance, of a
residual norm reduction to $\tol = 10^{-6}$ as the stopping criterion
for each time step, and experiment with the more practical algorithmic
alternative in Algorithm~\ref{alg:cgmres2}.
\revised{Figure~\ref{fig:lkdvevolve} presents two experiments.  At left, we consider the time evolution using a zero initial guess for each time step.  At right, we use the solution from the previous time step as the initial guess.  As shown in~\cite[Section 3.3]{BirkenLinders:2021}, using the solution from the previous time step leads to automatic conservation of mass, since we use the unpreconditioned algorithms here.  Thus, when we impose the constraints on CGMRES in this case, we only impose constraints on momentum and energy, since the conservation of mass constraint is already satisfied.}  
We observe in Figure~\ref{fig:lkdvevolve}
that all invariants are preserved to near machine
precision by CGMRES.
\revised{
For standard FGMRES, on the other hand, we see
that while mass is constrained to a relatively high tolerance in both cases, both
momentum and energy deviate up to the
order of the solver tolerance.}

\begin{figure}[h!]
  \caption{The error propagation over time using direct linear
    solvers, FGMRES and CGMRES using an $s$ stage GLRK temporal
    discretisation and the order $q$ spatial discretisation
    \eqref{eqn:lkdvrk} simulating \eqref{eqn:lkdv:exact}. We fix
    $\dt{}=0.1$, $\M=400$, and vary the solver tolerance $\tol$ to be
    $10^{-3}$, $10^{-5}$ and $10^{-7}$, respectively, as we increase the order of the
    method. Further, we precondition by \eqref{eqn:lkdv:pre} and
    initialise the linear solvers with either the stage values on the
    previous step, or a tiling of the initial data (for the first step).}
  \label{fig:lkdvHO}
  \includegraphics[width=0.9\textwidth]{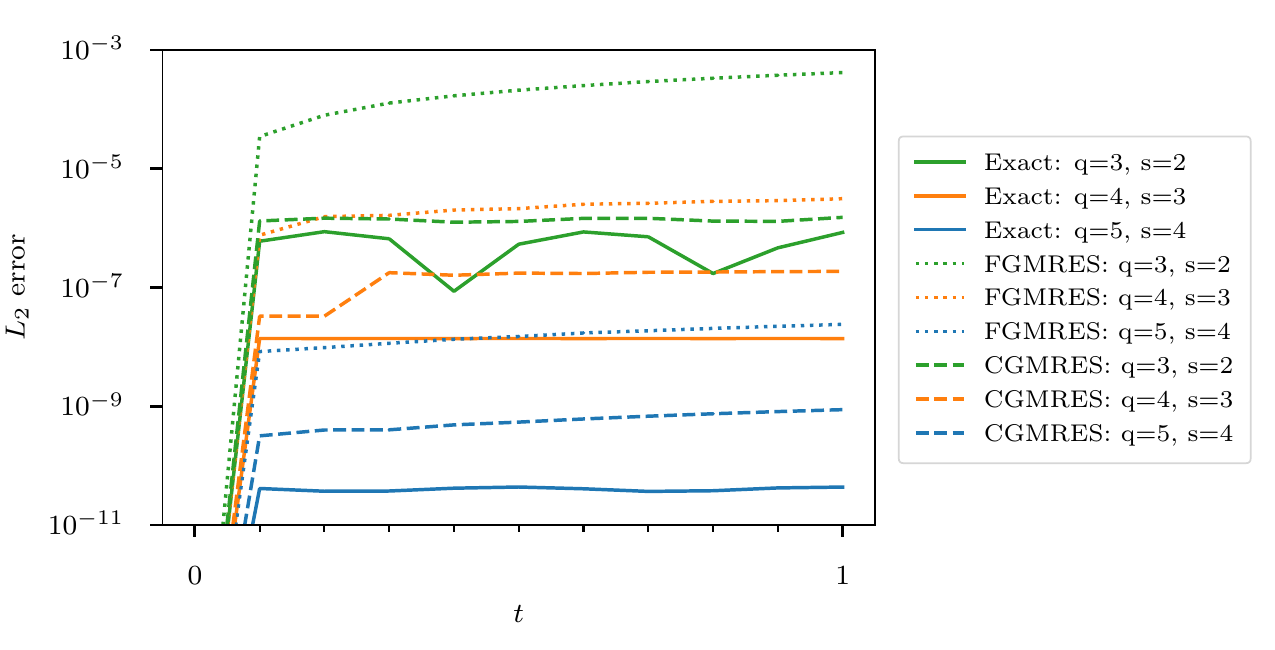}
\end{figure}

As discussed in Section~\ref{sec:background:lkdv}, the methodology proposed
here can be used for both \CN and the higher-order GLRK
discretisations given in
Definition~\ref{def:lkdvrk}.  Next, we consider performance as we
increase the temporal and spatial order, for a problem posed on a
larger domain and discretised with more points in space.  As is typical with any discretised PDE, increasing the
order of and number of points in the discretisation also increases the
condition number of the linear system to be solved, to the point where
effective convergence of GMRES-like methods requires the use of a
preconditioner.  Here, since we discretise a one-dimensional PDE, we
can make use of a simple preconditioner, such as the incomplete LU
(ILU) factorisation of the system.  Here, we make use of the ILU
implementation provided by SuperLU~\cite{superlu_ug99, lishao10},
through its interface to SciPy~\cite{2020SciPy-NMeth}, using a
supernodal ILUTP algorithm with drop tolerance $10^{-4}$ and fill
ratio upper bound of $10$, given by the call
\begin{equation} \label{eqn:lkdv:pre}
  \precon
  =
  \verb|scipy.sparse.linalg.spilu|(A,
  \verb|drop_tol| = 10^{-4},
  \verb|fill_factor| = 10)
  .
\end{equation}
Simulating the trigonometric wave
\begin{equation} \label{eqn:lkdv:exact}
  u(t,x)
  =
  \sin{\alpha\bc{x - (1-\alpha^2)t}} + 1
  ,
\end{equation}
with $\alpha=\frac{\pi}{5}$ over $t\in [0,1], x \in [0,40)$, for FGMRES,
CGMRES, and exact linear solvers, we obtain Figure~\ref{fig:lkdvHO},
showing the $L_2$ error in the numerical solution at each time step,
with $\dt{} = 0.1$.  We note here that we have greatly increased both
the size of the spatial domain and the number of points used to
discretise it, in comparison to the previous example.  Note, however,
that the spatial $L_2$ error in the discretisation cannot scale better
than $(X/\M{})^{q+1}$ (as degree $q$ polynomials are used)
while the temporal error for an $s$-stage GLRK
method should scale like $\tau^{2s}$, so we use slightly higher order
spatial discretisations than temporal ones.
For each discretisation, we plot the error observed using an exact
solver for each timestep (LU factorisation), preconditioned FGMRES,
and preconditioned CGMRES. \revised{We observe that CGMRES
  consistently outperforms FGMRES in this simulation, but remains less
  accurate than the exact solver.}  In fact, the solution error using
FGMRES can be larger than that of solving the next lower-order
discretisation using the exact solver, while CGMRES generally yields errors
between those from the exact solver at the same order and at one order lower.

\subsection{Shallow water equations} \label{sec:examples:swe}

\begin{figure}[h!]
  \caption{A comparison between CGMRES as \revised{proposed in Section~\ref{sec:pcgmres}} and FGMRES for the linear shallow water
    equations. Here, the linear system corresponds to the numerical
    approximation \eqref{eqn:swe:scheme} with ${\dt{}}=0.1$, ${\M}=50$, and
    variable degree $q$ and the constraints are described by
    \eqref{eqn:swe:constraints}. We initialise $\vec{U}^0$ and $\Rho^0$
    with
    \eqref{eqn:swe:ic} and solve for $(\vec{U}^1,\Rho^1)$. Further, we
    specify the initial guess for our linear solvers to be
    $\vec{x}_0=\vec{0}$.
    \label{fig:swe}}
  \includegraphics[width=0.9\textwidth]{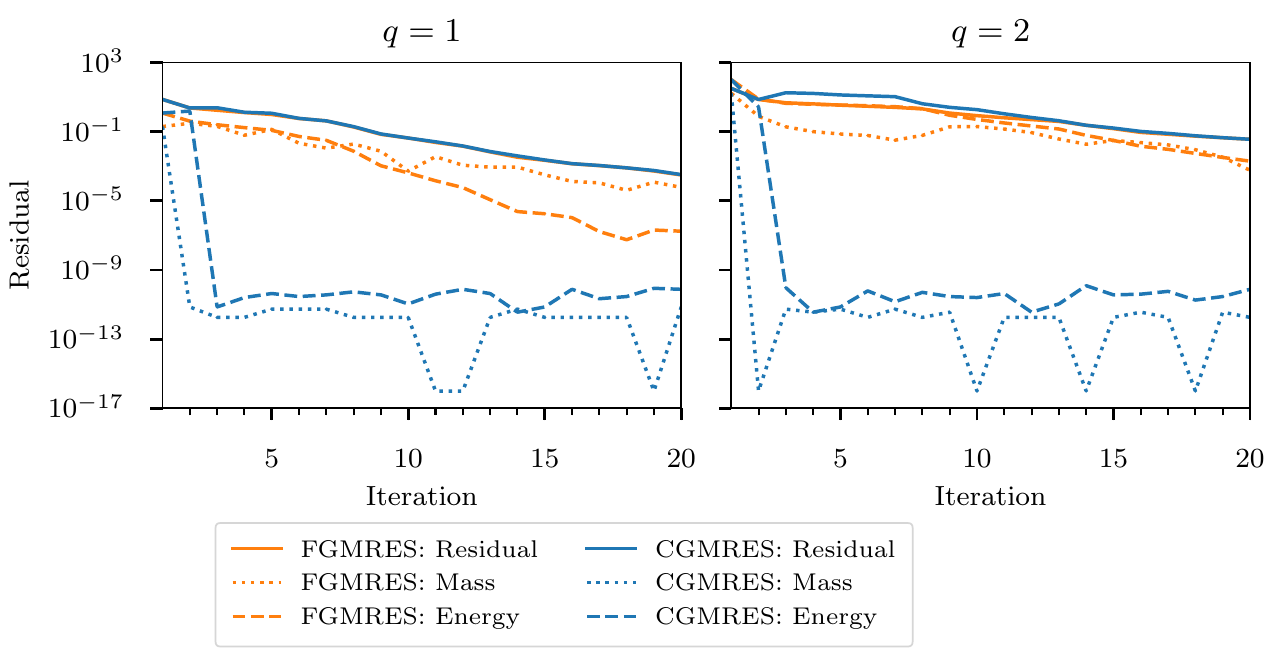}
\end{figure}

Here, we consider~\eqref{eq:swe1} with $c=1$ and $f=0.1$ over the domain $t\in[0,10], \vec{x} \in
[0,40) \times [0,40)$, where space is doubly periodic.
 We initialise our simulations with zero velocity and a Gaussian
 pressure distribution,
\begin{equation} \label{eqn:swe:ic}
  \begin{split}
    \vec{U}^0 & = \vec{0} \\
    \Rho^0 & = 10\, \exp{- \frac{\bc{x-20}^2 + \bc{y-20}^2}{20^2}}
    .
  \end{split}
\end{equation}
In Figure~\ref{fig:swe}, we again compare CGMRES against FGMRES in terms of
residual norm and the preservation of conserved quantities, for degree
$q=1$ and $q=2$, over 20 iterations. We note that the conserved quantities
converge at the same rate as the residual with FGMRES, whereas we are
able to
preserve both invariants at the level of machine precision for all
iterations after the third using CGMRES.
As for linear KdV, for early iterations, we see that the residual for
CGMRES is slightly larger than that for FGMRES, but that the two
residuals quickly become comparable.
For higher
polynomial degree the convergence of both iterations is significantly
slower, as is typical.  While we do not consider this in Figure~\ref{fig:swe},
convergence could be improved by either improving the initial guess
(we use a zero initial guess here, for illustration), or by the use of
a suitable preconditioner.

\begin{figure}[h!]
  \caption{The deviation in mass and energy over time for the finite
    element scheme for the shallow water equations
    \eqref{eqn:swe:scheme}, with ${\dt{}}=0.1$, $T=10$, ${\M}=50$, and $q=1$, using the linear solvers FGMRES and CGMRES as
    described by Algorithms \ref{alg:gmres} and \ref{alg:cgmres2},
    respectively. We constrain CGMRES by
    \eqref{eqn:swe:constraints} and set $\tol=10^{-6}$ and
    $\Epsilon=10\tol$. We initialise $\vec{U}^0$ and
    $\Rho^0$ with \eqref{eqn:swe:ic} and solve iteratively for
    $(\vec{U}^{\N},\Rho^{\N})$. \revised{At left, we present results for both constraints using a zero initial guess for the linear solvers at each time step.  At right, we present results for both constraints using the solution from the previous time step as the initial guess at each time step, but only enforcing the constraint on energy in CGMRES.}}
   \label{fig:sweevolve}
   \includegraphics[width=0.9\textwidth]{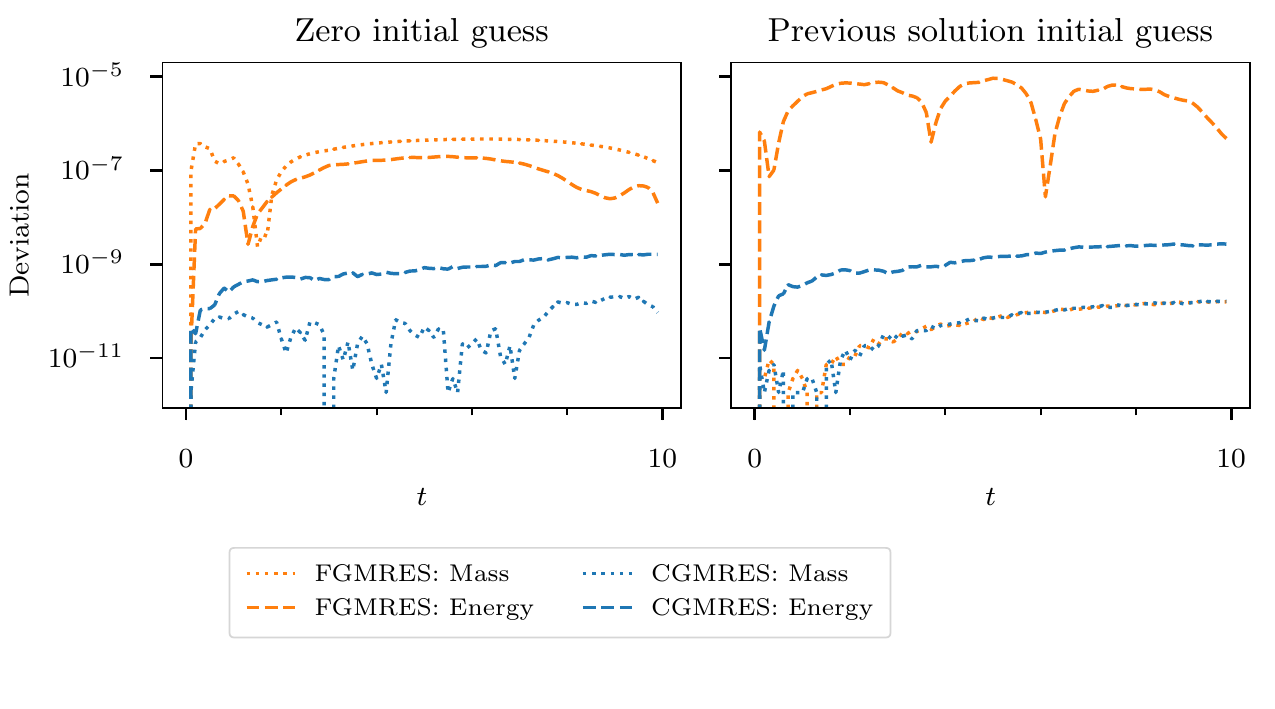}
\end{figure}

As for linear KdV, we next study the effect of
the different linear solvers for many timesteps, using a typical GMRES
stopping tolerance of
$\tol=10^{-6}$.  We consider the same problem as above, with $\dt{} =
0.1$ and $q=1$.
\revised{
As above, Figure~\ref{fig:sweevolve} presents results for both a zero initial guess and using the solution from the previous time step as the initial guess at each time step, showing that mass is conserved equally well in the latter case.
  }
In Figure~\ref{fig:sweevolve}, we observe that CGMRES
accurately preserves the conservation laws well below the linear
solver tolerance over long time, while energy conservation using FGMRES
appears limited to that given by the solver tolerance\revised{, as is mass conservation when using a zero initial guess}.  Here, we do
not have an analytical solution for comparison of accuracy between the
two solution schemes.

\revised{%
  Table~\ref{tab:swe} presents computational times for both FGMRES and CGMRES, using the same ILU preconditioner as considered above, but with drop tolerance set to $10^{-2}$ instead of $10^{-4}$.  We note that ILU is a rather expensive, but effective, preconditioner for this case, with the preconditioner setup time being, by far, the dominant time of the computation.  For these parameters, we see (as expected) some growth in the number of iterations as $\M$ increases, leading to faster-than-linear growth in time-to-solution with problem size (which is proportional to $\M^2$).  For each iteration of CGMRES, we time the cost per iteration, measuring the sum of the preconditioner application, matrix-vector product, Arnoldi orthogonalisation and solution of the Hessenberg system for the minimiser over the Krylov space.  This is reported as $\Titer$, the average time for iterations in which constraints are not imposed, and $\Titercon$, the average time for iterations in which the constraints are imposed.  Additionally, we measure the cost of projecting the constraints onto the Krylov space (see Remark~\ref{rem:cost}), again averaged over all constrained iterations, reported as $\Toverhead$ in Table~\ref{tab:swe}.  For the largest values of $\M$ reported in the table, we see that the difference between the total FGMRES and CGMRES solve times is almost exactly equal to the cost of assembling the constraints, which is bounded by about twice the average cost of a single iteration of CGMRES.  We note that our implementation of FGMRES does not implement the usual updating-QR factorisation approach to solving the Hessenberg system but, instead, uses the same} \verb|scipy.optimize.minimize| \revised{command to solve the system as in the constrained case, just without imposing the constraints.  We have compared timings to those of a more standard FGMRES, as implemented in} \verb|PyAMG|\revised{~\cite{BeOlSc2022} (noting that} \verb|SciPy| \revised{does not provide an FGMRES function), and found that our approach is slightly faster, at least for the small Krylov spaces we have used for comparison.  We also note that, as expected, the convergence of GMRES preconditioned by ILU suffers as the problem size gets larger and the convergence tolerance gets stricter.  Thus, while the added cost of CGMRES over FGMRES is noticeable in these examples (20-30\% for $\M = 256$ and $512$), this cost can clearly pay off when ``regular'' FGMRES requires many more iterations to enforce the constraints to suitable tolerances than are needed for the expected/desired residual reduction.
}

\begin{table}[h!]
  \centering
  \revised{
    \caption{\revised{%
        Computational times for a single solve of the
        linear shallow water equations \eqref{eqn:swe:scheme}
        preconditioned by \eqref{eqn:lkdv:pre} for variable spatial
        mesh resolutions $\M$, with $\dt{}=0.1$, $\tol=10^{-7}$, $\Epsilon=10\tol$ and
        a zero initial guess for $\vec{x}_0$. We present the preconditioner assembly
        time (in seconds) $\Tprecon$, total run time of FGMRES (as
        described in Algorithm \ref{alg:gmres}) $\Tgmres$, total
        run time of CGMRES (as described in Algorithm
        \ref{alg:cgmres2}) $\Tcgmres$, and total number of iterations for both FGMRES and CGMRES (which are the same in this example). We further break down the costs within
        CGMRES, writing $\Titer$ as the average run time for an
        unconstrained iteration, $\Toverhead$ as the average
        computational overhead (per constrained iteration) required to
        assemble the constraints, and $\Titercon$ as the constrained
        iteration cost. Here, $\len{\Titercon}$ denotes the number of constrained iterations of CGMRES.
        \label{tab:swe}}}
    \begin{tabular}{||c||c|c|c|c|c||}
      \hhline{#=#=|=|=|=|=#}
      $\M$ & 32 & 64 & 128 & 256 & 512 \\
      \hhline{#=#=|=|=|=|=#}
      $\Tprecon$ & 4.18e-02 & 1.63e-01 & 8.96e-01 & 5.26e+00 & 3.07e+01 \\
      \hline
      $\Tgmres$ & 1.18e-02 &4.42e-02 &  2.08e-01 &1.09e+00 & 7.03e+00 \\
      \hline
      $\Tcgmres$ & 2.24e-02 & 7.90e-02 &5.50e-01 & 1.38e+00 &  8.48e+00 \\
      \hline
      Iterations & 6 & 6 & 6 & 7 & 9 \\
      \hhline{#=#=|=|=|=|=#}
      $\Titer$ &  2.31e-03 & 7.43e-03 & 3.40e-02 &  1.51e-01 & 7.64e-01 \\
      \hline
      $\Toverhead$ &3.53e-03 &  1.25e-02 &  4.84e-02 & 2.90e-01 &  1.43e+00 \\
      \hline
      $\Titercon$ & 6.95e-03 &  2.70e-02 & 3.27e-01 & 1.63e-01 &  8.41e-01 \\
      \hline
      $\len{\Titercon}$ & 1 & 1 & 1 & 1 & 1 \\
      \hhline{#=#=|=|=|=|=#}
    \end{tabular}
  }
\end{table}

\subsection{Heat equation} \label{sec:examples:heat}

Recall the heat equation and associated discretisation discussed
in Section~\ref{sec:background:heat}. 
Experimentally, we consider a single time step with initial condition
given by
\begin{equation} \label{eqn:heat:ic}
  U^0
  =
  \Pi \bc{
    10^3 \bs{ \bc{ x (x-1)}^5 + \bc{y (y-1)^6}}
  }
  ,
\end{equation}
over the domain $t \in [0,1]$, $\vec{x} \in [0,1]\times[0,1]$,
\revised{where $\Pi$ is the $L_2$ projection into the finite element
  space. Spatially, we employ Neumann boundary conditions.} Due to the
poor conditioning of the problem, we expect unpreconditioned iterative
solvers to be slow to converge. With this in mind, \revised{we seek to use an effective preconditioner in the iteration.
  Unsurprisingly, performance using the ILU preconditioner deteriorates significantly when we increase problem size.  Thus, for a timing and scaling study, we consider a much more robust preconditioning framework for the heat equation, that of Ruge-St\"{u}ben algebraic multigrid~\cite{ABrandt_SFMcCormick_JWRuge_1984a,Ruge1987,KStuben_2001a}.

  Multigrid
  methods~\cite{WLBriggs_VEHenson_SFMcCormick_2000a,UTrottenberg_etal_2001a}
  are widely recognised as among the most efficient families of
  preconditioners for elliptic PDEs and for simple time
  discretisations of parabolic PDEs (although they can be readily
  adapted for more complex problems and discretisations;
  cf.~\cite{Irksome:2021,RAbuLabdeh_etal_2022a}), combining fine-scale
  relaxation schemes with coarse-level correction.  For many problems,
  \emph{geometric} multigrid methods, where the multigrid hierarchy is
  created using uniform refinement of structured grids with geometric
  (e.g., piecewise linear) grid-transfer operators, are the most
  efficient approaches.  However, \emph{algebraic} multigrid methods
  offer black-box solvers that can be easily applied based on the
  system matrix alone, without appealing to any geometric information
  about the underlying discretisation or meshes.  Here, we use the
  implementation of Ruge-St\"{u}ben algebraic multigrid (AMG) with
  default parameters from} \verb|PyAMG|\revised{~\cite{BeOlSc2022},
  determining coarse grids using the first pass of the original
  Ruge-St\"uben coarsening algorithm, with standard
  strength-of-connection using parameter $\theta = 0.25$, classical
  interpolation, and a V(1,1) cycle with symmetric Gauss-Seidel as the
  pre- and post-relaxation scheme.  }

\begin{figure}[h!]
 \centering
 \caption{ A comparison between CGMRES as \revised{proposed in Section~\ref{sec:pcgmres}} and FGMRES for the heat equation. The linear system corresponds to the numerical scheme
   \eqref{eqn:heat:scheme} with ${\dt{}}=0.01$, ${\M}=50$ and degree
   $q=1$. We initialise $U^0$ with \eqref{eqn:heat:ic} and solve for
   $U^1$, enforcing the constraints \eqref{eqn:heat:constraints} and
   using initial guess $\vec{x}_0 = \vec{0}$. Here, we study the
   effects of preconditioning the linear system \revised{using algebraic multigrid (shown at right).}}
 \label{fig:heat}

 \includegraphics[width=0.9\textwidth]{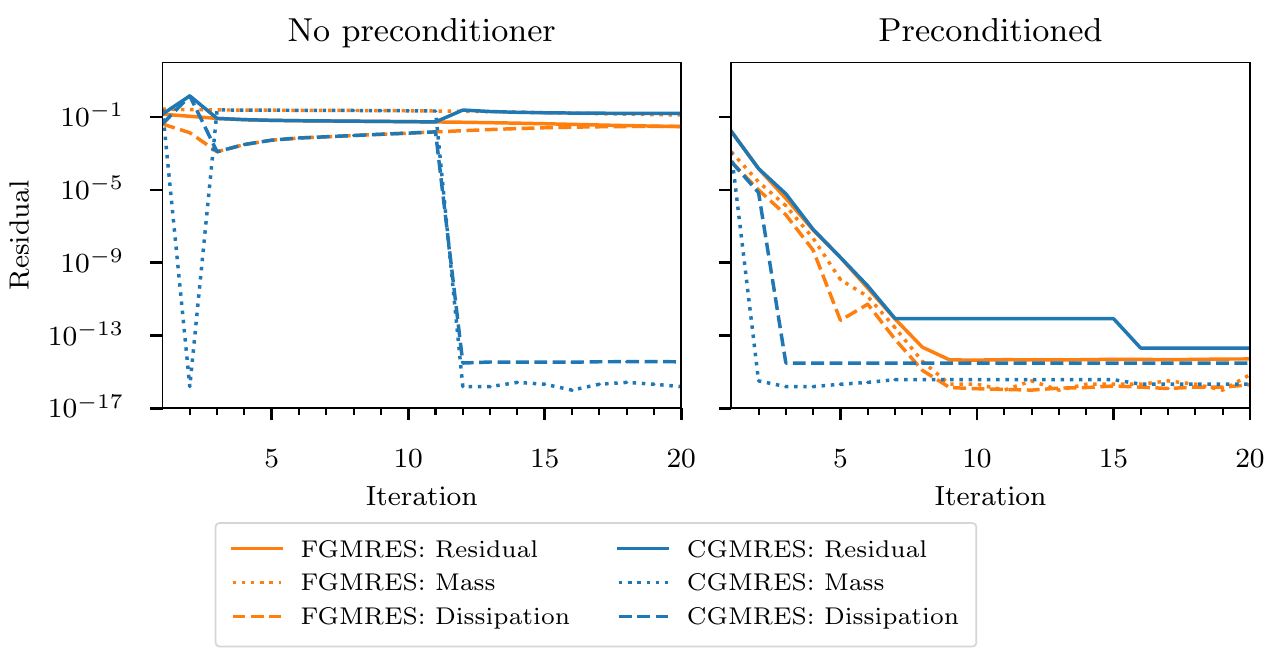}
 
\end{figure}

In Figure~\ref{fig:heat}, we compare CGMRES against FGMRES for both the
standard and preconditioned linear system, for a fixed number of iterations. In the nonpreconditioned
case, both algorithms stagnate in residual.  We see that enforcing the
dissipation law on CGMRES leads both to an increase in residual norm over
the first tens of iterations (in comparison with that of FGMRES),
while failing to satisfy the constraints until iteration $12$.
\revised{
  Here, we note that the }\verb|SQSLP| \revised{solver returns not-a-numbers when it fails to converge, leading to our inclusion of the check on Line~\ref{alg:cgmres2:check} of Algorithm~\ref{alg:cgmres2}.  When this happens, we revert to the regular FGMRES solution computed on Line~\ref{alg:cgmres2:redirect}, explaining the concurrence between the residual norms for iterations 3-12.  Note also that, for these iterations, the misfit in the mass conservation is larger than the residual norm, showing there is no general relationship between the size of the residual norm and the misfit in the constraints.
  }
Even
once the Krylov space is rich enough to allow the constraints to be
satisfied, we are still far from convergence of the linear system.  In
contrast, we see similar reduction in residual norm for both CGMRES
and FGMRES with the \revised{AMG} preconditioner.  Furthermore, the CGMRES
iteration shows immediate satisfaction of the constraints when they
are imposed, with no significant deviation in the residual norm until
we are well past reasonable stopping tolerances, with the residual norm
below $10^{-9}$.

\revised{%
  Results in~\cref{tab:heat} again show the relative costs of preconditioner setup and solution using classical FGMRES and the CGMRES algorithm proposed here.  In particular, we now note that the AMG setup time is on the same order as both solve times (as is typical for AMG), and that the additional cost for the CGMRES iterations is again primarily due to the overhead of computing the ``reduced'' constraints over the Krylov space.  This cost scales, as expected, with problem size ($\M^2$) and size of the Krylov space, and is roughly equal to the average cost of a single step of the unconstrained algorithm.  Here, we note that the preservation of the constraints is significantly better with CGMRES than with FGMRES, particularly for the mass constraint.  When running FGMRES, the constraints generally converge roughly at the same rate as the GMRES residual norm.  As seen in Figure~\ref{fig:heat}, this can lead to a large difference in the misfit in the constraints between FGMRES and CGMRES when we converge in 5-6 iterations, as we do here.  Thus, while it is true that we could substantially reduce the FGMRES constraint misfit by running one or two more iterations of FGMRES, the cost of doing so would equal or exceed that of the extra overhead for enforcing the constraints once at the end of a CGMRES solve.
}

\begin{table}[h!]
  \centering
  \revised{
    \caption{\revised{%
        Computational times for a single solve of the heat equation
        \eqref{eqn:heat:scheme} preconditioned by Ruge-Stuben AMG
        for variable spatial mesh resolution, $\M$, with $\dt{}=0.1$, $\tol=10^{-7}$, $\Epsilon=10\tol$ and
        a zero initial guess for $\vec{x}_0$. We present the preconditioner assembly
        time (in seconds) $\Tprecon$, total run time of FGMRES (as
        described in Algorithm \ref{alg:gmres}) $\Tgmres$, total
        run time of CGMRES (as described in Algorithm
        \ref{alg:cgmres2}) $\Tcgmres$, and total number of iterations for both FGMRES and CGMRES (which are the same in this example). We further break down the costs within
        CGMRES, writing $\Titer$ as the average run time for an
        unconstrained iteration, $\Toverhead$ as the average
        computational overhead (per constrained iteration) required to
        assemble the constraints, and $\Titercon$ as the constrained
        iteration cost. Here, $\len{\Titercon}$ denotes the number of constrained iterations of CGMRES.
        \label{tab:heat}}}
    \begin{tabular}{||c||c|c|c|c|c||}
      \hhline{#=#=|=|=|=|=#}
      $\M$ & 128 & 256 & 512 & 1024 & 2048 \\
      \hhline{#=#=|=|=|=|=#}
      $\Tprecon$ &  1.93e-02 &  5.81e-02 &   2.35e-01 & 1.15e+00 &  4.54e+00 \\
      \hline
      $\Tgmres$ &  1.69e-02 & 5.86e-02 &  2.62e-01 &  1.19e+00 &  6.50e+00 \\
      \hline
      $\Tcgmres$ &  2.28e-02 &  7.13e-02 &  3.07e-01 &   1.41e+00 &  8.60e+00 \\
      \hline
      Iterations & 5 & 5 & 5 & 5 & 6 \\
      \hhline{#=#=|=|=|=|=#}
      $\Titer$ &   3.72e-03 &  1.18e-02 &  5.12e-02 &  2.30e-01 &  1.04e+00 \\
      \hline
      $\Toverhead$ & 3.04e-03 & 1.03e-02 &   4.30e-02 &  2.35e-01 &  1.05e+00 \\
      \hline
      $\Titercon$ &  4.50e-03 &  1.28e-02 &  5.49e-02 &   2.41e-01 &  1.12e+00 \\
      \hline
      $\len{\Titercon}$ & 1 & 1 & 1 & 1 & 2 \\
      \hhline{#=#=|=|=|=|=#}
    \end{tabular}
  }
\end{table}

\FloatBarrier
\section{Conclusion}\label{sec:conclusion}

Here, we introduce a modification to the FGMRES algorithm that allows
for conserved quantities of an
underlying discretisation to be preserved at any desired stopping
tolerance. While we focus on the modification for FGMRES, the approach can be freely applied to any Krylov solver
where an explicit basis for the Krylov space is used in a
minimisation algorithm. Importantly, the proposed
constrained algorithm can be incorporated into existing solver
implementations through minor modification, replacing a QR
factorisation of the $(\ell+1)\times \ell$ Hessenberg matrix with a
(nonlinear) constrained optimisation problem in $\ell+c$ unknowns for
$c$ constraints. We observe
(experimentally) that the constrained solver performs no worse than
FGMRES, so long as the Krylov space is sufficiently rich and/or the
initial guess is sufficiently close to solution.  While we focus here
on evaluating the method as a proof-of-concept, we believe the experiments are promising enough to
justify future work on efficient implementation of this approach in an
optimised high-performance computing environment.  Such a study could
also include the nonlinear variant of this approach, proposed in Remark~\ref{rem:nonlinear}.

\FloatBarrier
\bibliographystyle{siamplain}
\bibliography{geosolve}

\end{document}